 \documentclass{article}
 \usepackage{times}
 \usepackage{fullpage}
 \usepackage{geometry}
 \usepackage{amsmath}
 \usepackage{amsthm}
 \theoremstyle{definition}

 \newcommand{\EquationsNumberedThrough}{}
 \newcommand{\TheoremsNumberedThrough}{}
 \newcommand{\ECRepeatTheorems}{}
 \newcommand{\MANUSCRIPTNO}[1]{}
 \newcommand{\RUNAUTHOR}[1]{}
 \newcommand{\RUNTITLE}[1]{}
 \newcommand{\TITLE}[1]{\title{#1}}
 \newcommand{\ARTICLEAUTHORS}[1]{\author{#1}}
 \newcommand{\AUTHOR}[1]{#1}
 \newcommand{\AFF}[1]{\thanks{\small{#1}\vspace{.66em}} \and}
 \newcommand{\EMAIL}[1]{\\\small{\texttt{#1}}}
 \newcommand{\ABSTRACT}[1]{\begin{abstract}\noindent #1\end{abstract}}
 \newcommand{\FUNDING}[1]{\noindent\textbf{Funding.} #1\\}
 \newcommand{\KEYWORDS}[1]{\noindent\textbf{Keywords.} #1\\}

 \newcommand{\APPENDIX}{\appendix}

 \newcommand{\TABLE}[3]{#2 \caption{#1}}
  \newcommand{\FIGURE}[3]{#1 \caption{#2}}
 
\usepackage[switch]{lineno}

\newcommand{\ignore}[1]{}

\usepackage{endnotes}
\let\footnote=\endnote
\let\enotesize=\normalsize
\def\notesname{Endnotes}%
\def\makeenmark{$^{\theenmark}$}
\def\enoteformat{\rightskip0pt\leftskip0pt\parindent=1.75em
	\leavevmode\llap{\theenmark.\enskip}}

\usepackage{graphicx}
\usepackage{eqndefns-left}
\usepackage{multirow}
\usepackage{hhline}
\usepackage{longtable}
\usepackage{placeins}
\usepackage{booktabs}
\usepackage{pgfplots}

\usepackage[small, margin=1cm]{caption}

\usepackage{appendix}
\usepackage{color}
\definecolor{strcolor}{rgb}{0.6, 0.2, 0.6}
\definecolor{commentcolor}{rgb}{0.3125, 0.5, 0.3125}
\definecolor{keycol}{rgb}{0, 0, 1}

\usepackage{amsmath}

\usepackage{listings}
\usepackage{url}

\def\blot{\quad \mbox{$\vcenter{ \vbox{ \hrule height.4pt
				\hbox{\vrule width.4pt height.9ex \kern.9ex \vrule width.4pt}
				\hrule height.4pt}}$}}

\usepackage{natbib}
\bibpunct[, ]{(}{)}{,}{a}{}{,}%
\TheoremsNumberedThrough     
\ECRepeatTheorems

\EquationsNumberedThrough    

\gdef\AQ#1{}
\gdef\CQ#1{}

\begin{document}

	\RUNAUTHOR{Blom et~al.} %

	\RUNTITLE{Long-Term Open-Pit Mine Planning with Large Neighbourhood Search}

\TITLE{Long-Term Open-Pit Mine Planning with Large Neighbourhood Search}

	\ARTICLEAUTHORS{
\AUTHOR{Michelle Blom, Adrian R. Pearce}
\AFF{School of Computing and Information Systems, The University of Melbourne, Parkville, Victoria, Australia \EMAIL{michelle.blom@unimelb.edu.au, adrianrp@unimelb.edu.au}}
\AUTHOR{Pascal C\^{o}t\'{e}}
\AFF{Operations Research and Decision Science, Rio Tinto, Canada \EMAIL{Pascal.Cote@riotinto.com}}
}

\date{}

\maketitle

\ABSTRACT{
We present a Large Neighbourhood Search  based approach for solving complex long-term open-pit mine planning problems. An initial feasible solution, generated by a sliding windows heuristic, is improved through repeated solves of a restricted mixed-integer program. Each iteration leaves only a subset of the variables in the  planning model free to take on new values. We form these subsets through the use of neighbourhood formation strategies that exploit model  structure. We show that our approach is able to find near-optimal solutions to problems that cannot be solved by an off-the-shelf solver in a reasonable time frame, or with reasonable computational resources. Our method substantially reduces the solve times required for large models, allowing mine planners to explore multiple scenarios in a timely fashion. Our approach is being used by Rio Tinto to solve large long-term mine planning problems, and has been responsible for generating millions of dollars in value insights.
}

\FUNDING{This work was
supported by Rio Tinto, the University of Melbourne's Research Computing Services, and by the Australian Research Council (OPTIMA ITTC IC200100009).}


Rio Tinto is a British-Australian multinational mining and materials cooperation with headquarters in London, United Kingdom, and Melbourne, Australia.   Formed in 1873 with the purchase of the Rio Tinto mine complex in Spain, Rio Tinto is now the second largest mining and materials company in the world. The company operates in 35 countries, with a portfolio that includes iron ore, aluminum, copper, borates, diamonds, lithium, salt, scandium, and titanium dioxide. Rio Tinto is one of the worlds leading producers and exporters of iron ore through its operations in the Pilbara region of Western Australia. Since the opening of Mount Tom Price in 1966, iron ore production in the Pilbara has, in 2024,  grown to over 320 million tonnes per year through a network of 17 mines connected by 2000 km of rail to 4 port terminals. Rio Tinto produces several iron ore products, including the high grade Pilbara blend, formed by blending ore from mines across this network.

\section{Problem Summary}

Rio Tinto has invested in the development of in-house software to optimise long-term open-pit mine plans. This software translates an open-pit block model and supply chain network design into a mixed-integer linear programming model. For small to medium-sized problems, commercial solvers are able to solve the resulting models in reasonable time frames (an hour or less). The use of the in-house platform on larger deposits and more complex mining systems resulted in a substantial increase in model sizes, with solve times becoming a major concern. In practice, mining engineers using this platform solve a stochastic problem by solving many instances of a deterministic one. Across these instances, parameters that capture aspects such as price and grade are varied to reflect the uncertainty present in the original problem.  Key strategic decisions are made by analysing the resulting plans across these scenarios.

Mining engineers  began simplifying the formulation of large models by removing features such as stockpiles and truck fleets. However, this produced mine plans that underestimated  production capacity and flexibility. Decomposition strategies were considered, such as removing stockpiles, solving the problem, fixing part of the block sequence, reintroducing stockpiles, and then performing a final solve. These strategies produced unrealistic mine plans or infeasible solutions. For the models considered in this paper, the presence of stockpiles adds substantial complexity. These stockpiles are required, however, to satisfy minimum production constraints across the planning horizon. For large models not amenable to full solves with a commercial solver,  mining engineers apply heuristics like the sliding windows method. This technique iteratively builds a complete mine plan while solving the problem for  a small number of periods in a window that slides along the horizon after each solve.    

With the extension of the platform to large Rio Tinto systems like the Pilbara, the need to develop a faster solver that is able to find near-optimal long-term plans has become critical for the software's acceptance by mining engineers. This has provided motivation for alternative approaches.

This paper considers a long-term production scheduling problem for a set of open-pit mines, connected by rail to junctions at which blending takes place. 
Material from different regions of a mine site (called `blocks’) is extracted, in each period of a horizon, and sent to one of a number of destinations (e.g., processing plants, waste dumps, and stockpiles). Precedences constrain the order in which blocks can be extracted and capacity constraints limit the total tons of material mined in each period, and sent to each destination. Blending takes ore of varying characteristics, such as grade, and levels of contaminants, and mixes it to produce a homogeneous product. Our model includes blending constraints, minimum  production constraints across regions of each mine (called pits), and capex decisions relating to the opening of pits.

We present a Large Neighbourhood Search based approach for solving long-term open-pit mine planning models, and a demonstration of the method on two real world case studies. The value of this method is that it allows us to find near-optimal solutions to large models without the need for simplifications that reduce realism and the subsequent value of solutions for decision-making.

 \section{Industrial Case Studies}\label{sec:Data}
We consider two case studies: the Oyu Tolgoi deposit in Mongolia; and the Pilbara mining system in Western Australia.

\subsection{Case study 1: Oyu Tolgoi}
\begin{figure}[t]
\centering
\FIGURE
{\includegraphics[width=\textwidth]{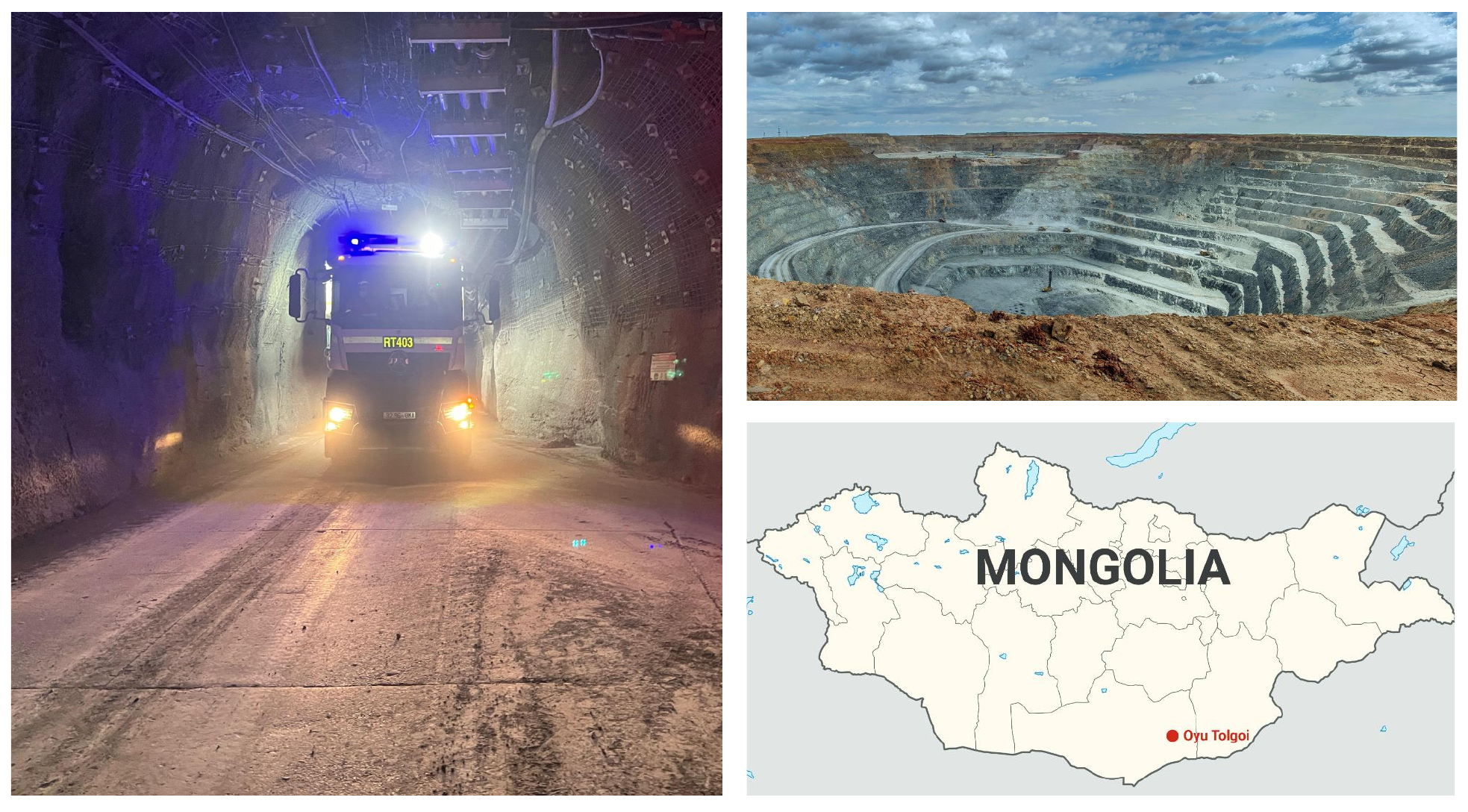}}
{Location, and operations, of the Oyu Tolgoi mining system. Image courtesy of Rio Tinto.\label{fig:OTOps}}
{}
 \end{figure}

 Oyu Tolgoi is a large copper and gold deposit in the south of the Gobi region in Mongolia, jointly owned by the Mongolian government and Rio Tinto. The deposit was discovered in the early 1990s by Canadian geologist John Brock. Robert Friedland's Ivanhoe Mines (later renamed Turquoise Hill Resources) began exploration in 2000, at which point the true scale of the deposit was revealed. Rio Tinto acquired a 66\% stake in the deposit through its purchase of Turquoise Hill Resources in December 2022, with the  Mongolian government holding the remaining 34\%. 
 The mining complex encompasses both open-pit and underground operations. Open-pit mining commenced in 2011, with the copper concentrator processing ore from 2013.
 The underground expansion, approved in January 2022, employs  block-caving techniques to access deeper ore bodies.
 At full capacity, Oyu Tolgoi is expected to produce approximately 500,000 tonnes of copper annually between 2028 and 2036. Figure \ref{fig:OTOps} shows the geographical location of  Oyu Tolgoi and images from the site.

This paper captures only the open-pit mines of the deposit. The model contains two large open pits, each composed of more than 100 blocks, and each block with parcels of multiple material types. The model has a supply chain network comprised of multiple dumps, stockpiles and crushers. Blending constraints ensure the production of saleable gold, silver, and copper products. Minimum production constraints are not present in this model. The Oyu Tolgoi model is defined over a horizon of 15 years, and contains 3.4 million variables, 2.6 million constraints, and a total of 9,377 parcels of material.

\subsection{Case study 2: Pilbara Iron Ore System}
Rio Tinto owns a portfolio of iron ore assets in the Pilbara region of Western Australia. In 1962, Rio Tinto discovered rich deposits at Mount Tom Price. Throughout the 1970s and 1980s, Rio Tinto expanded its operations significantly, opening new mines like Paraburdoo and Brockman. In the early 2000s, the global demand for steel surged, particularly due to China’s rapid industrialisation and urbanisation. This drove Rio Tinto to further expand its Pilbara operations. In the  2000's, Rio Tinto consolidated its position by merging and forming joint ventures with other players, including Robe River Iron Associates and Hope Downs. Rio Tinto operates an extensive rail network to transport iron ore from its mines to its ports on the Pilbara coast.

 \begin{figure}[t]
\centering
\FIGURE
{\includegraphics[width=\textwidth]{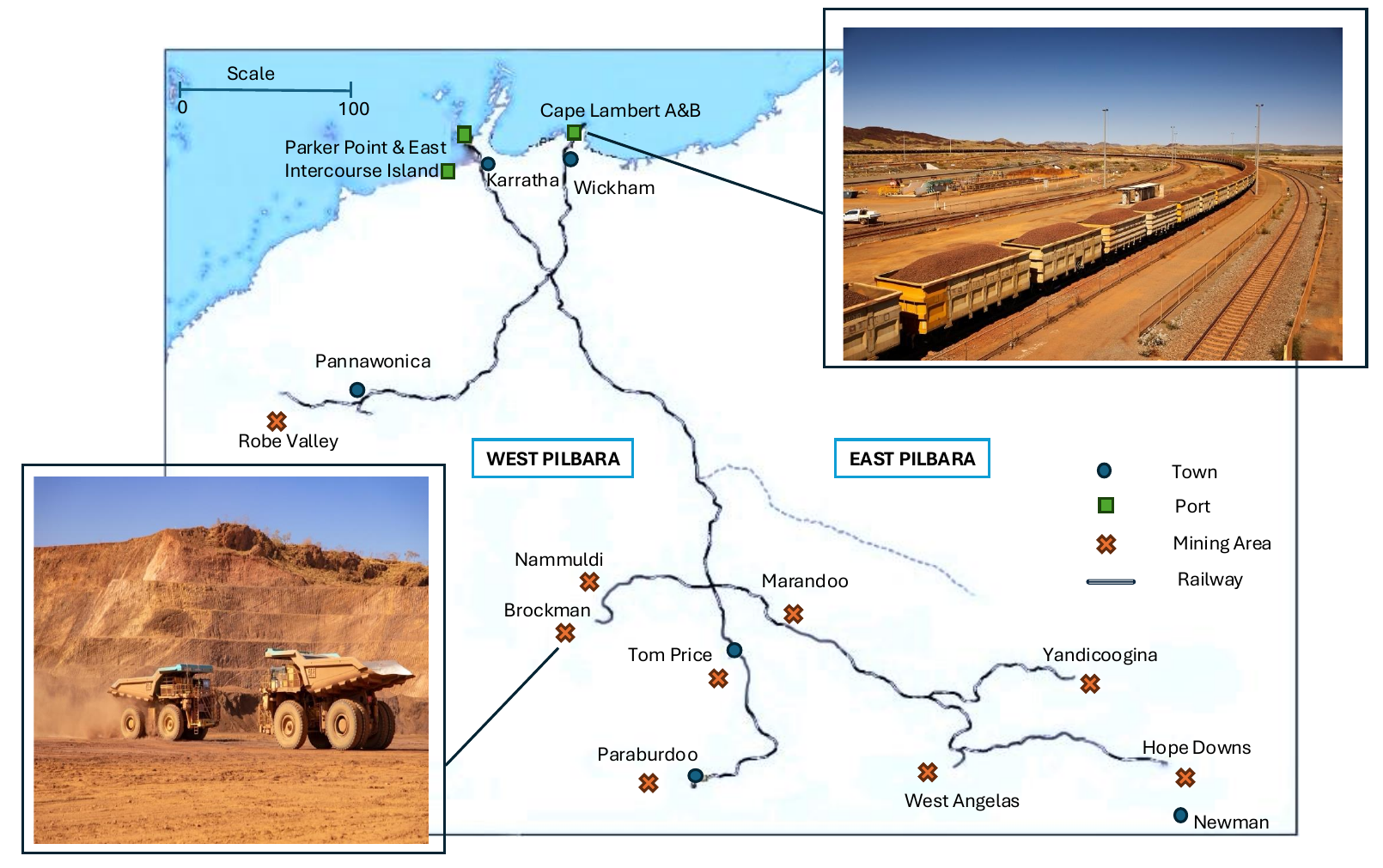}}
{The Pilbara mining system in Western Australia. Image courtesy of Rio Tinto.\label{fig:PilbaraOps}}
{}
 \end{figure}

The Pilbara iron ore system is composed of over 200 pits distributed across 17 mines, a selection of which are highlighted in Figure \ref{fig:PilbaraOps}. These pits contribute to five iron ore products whose composition is controlled by blending constraints. These constraints enforce a minimum concentration of iron, while limiting the concentration of alumina and phosphorus contaminants. These products are transported through a 2000 km rail system that connects the 17 hubs to 4 ports. Minimum  production constraints are enforced across subsets of pits. The model contains capex decisions relating to the timing of opening pits over a horizon of 15 years,  and contains 6.6 million variables, 4.1 million constraints, and  28,582 parcels of material.

\section{Solution Approach: Overview}

 Large Neighbourhood Search (LNS)  improves upon an existing solution to a problem  by repeatedly solving a \textit{simplification} in which a subset of decision variables are fixed to their values in the existing solution and the remainder left to take on a new value. The variables that are \textit{free} to take on a new value form a \textit{neighbourhood} around the current solution. The result is a smaller, restricted subproblem. Solving this subproblem allows us to search for a new and improved solution in this neighbourhood. 

Our approach first applies sliding windows  to find an initial feasible solution to a mixed integer linear programming model (MIP) of a long-term open-pit mine planning problem.  We then perform a \textit{parallelised} form of LNS over multiple threads of computation. On each thread, we form a neighbourhood with respect to the current best found solution. In our context,  this  neighbourhood is a collection of blocks. 
We leave the decision variables related to the mining of blocks in our neighbourhood free in our subproblem, and fix those decisions for blocks outside of the neighbourhood to their values in the current best solution. We  solve these restricted subproblems, updating our current best solution to the best found over all threads. 
We repeat these iterations until a termination criterion is satisfied. 

  To limit subproblem size, and to increase the likelihood that each neighbourhood will lead  to an improved solution, we use a path-based neighbourhood structure, and a series of neighbourhood formation strategies that exploit model structure. We form a neighbourhood by repeatedly selecting a block from the model, forming a path around that block, and adding the blocks in that path to the neighbourhood. We consider a range of strategies for selecting the blocks around which to form paths. These strategies define different methods for weighting the blocks in a model according to their value for scheduling purposes, and using those weights to guide selection. We outline how we can leverage what we know about different sets of model constraints to form neighbourhoods that increase the likelihood of finding new feasible solutions.  
We show that our approach substantially improves the quality of the initial solutions formed by sliding windows.
\section{Related Work}\label{sec:RelatedWork}

The scheduling of mines over long-term horizons (i.e., decades) has been well studied (see \cite{osa08}, \cite{new10}, \cite{ask11JOMS}, \cite{eps12}, \cite{lam14}, \cite{Lamghari2017-xk}, \cite{ZENG2021102274}, and \cite{Fathollahzadeh2021-jo} for reviews). \citet{eps12}  present MIP models characterising a variety of long-term scheduling problems. Of these, the precedence-constrained production scheduling problem (PCPSP) is most similar to the problem this paper tackles. 
In the PCPSP, blocks are extracted and sent to one of a number of facilities such as processing plants or stockpiles. Mining precedences constrain the order in which blocks are extracted, while side constraints can be specified to enforce resource capacities.

Many approaches have been developed to solve the PCPSP. 
\citet{cul11} present a sliding windows heuristic  for solving a variation of the problem in which a sequence of integer programs are solved. For each solve, all problem constraints are enforced in a window spanning the first few periods of the horizon. A Lagrangian relaxation of the model is enforced outside of this window. The resulting solution is used to fix the activities of the first period, after which the window slides forward by one period. The decision variables in subsequent solves are fixed for the first period, subject to all problem constraints in the periods covered by the window, and a relaxation of the model thereafter. The heuristic terminates once the last period in the horizon is scheduled.         
The sliding windows heuristic used in this paper is modelled on the work of \citet{cul11}, \citet{doi:10.1287/inte.2014.0737}, and \citet{o2015optimization}.

Much existing work applying neighbourhood-based techniques for mine planning uses simple neighbourhood structures. With an initial solution formed using the heuristic of \citet{gersh87}, 
\citet{Sari2016-bc} employ simple perturbations that change the period in which a selected block is mined, with a certain probability. 
\citet{Goodfellow2016-uk} employ three neighbourhood structures when solving a stochastic mine planning problem: randomly selecting a block to remove from the schedule or alter its time of extraction; altering the destination of material from a randomly selected cluster of blocks; and altering the subsequent destination of material leaving stockpiles or processing plants. \citet{LAMGHARI2020104590} present hyperheuristics that select among a range of simple perturbations: shift the mining of a block to the previous or subsequent time period; move the mining of block to another period; or swap the periods in which two blocks are mined. A range of variations on how blocks and time periods are selected are presented, forming 27 ways to perturb a schedule. 

\citet{ama09} solve a series of integer programs to incrementally improve an initial schedule found using the greedy heuristic of \citet{gersh87}. All model variables  are first fixed to their value in this initial schedule. A selection of variables are then unfixed and the resulting model solved to yield a (potentially) improved solution. These variables are then fixed to their value in this new solution. This process is repeated for all collections of variables generated in accordance with several strategies. The `cone above' strategy defines, for each block $b$, a set of variables relating to the mining of $b$ and its predecessors. The `periods' strategy defines, for each period $t$, a set of variables relating to the activities of $t$. The `transversal' strategy defines, for each block $b$, variables relating to the mining of all blocks within a defined distance of $b$. \citet{chic12} extend the heuristic of \citet{ama09} with different methods for selecting variables to fix. These methods select a random block $b$ that has been scheduled for mining in some period $t$, then: unfixes variables relating to the mining of $b$ and its predecessors (or, alternately, its successors); or unfixes all variables (to a maximum number) relating to blocks mined in periods $t-1$, $t$, and $t+1$. 

Our approach is similar to that of \citet{ama09}, in that we repeatedly select blocks around which to form a neighbourhood. Variables related to the mining of these blocks are left unfixed, and the resulting subproblem solved to form a new solution. Our approach differs in the type of neighbourhood structure used, and the strategies applied to instantiate it. Simple perturbative neighbourhood structures, such as swapping or shifting the periods in which blocks are mined do not form large enough neighbourhoods to allow new feasible solutions to be discovered for our model. Altering the mining of a block may require changes to how its predecessors and their predecessors, or its successors and their successors, are mined. This motivates the use of a path-based neighbourhood structure. This type of neighbourhood can be considered as a refinement of the often used cone-above or below structures. The strategies we use to select blocks around which to form paths focus on the constraints in our model that are tight, or most often responsible for infeasibilities (blending and minimum production constraints), and decisions that are likely to have a substantial impact on the objective (such as when pits are opened).

\citet{lam12} apply Tabu search to improve an initial schedule generated by greedily selecting eligible blocks to be mined in each time period. Their method  repeatedly shifts, adds, or deletes the mining of a block,  to or from an eligible period. Shifts that reverse recently performed actions form part of a Tabu-list, and are not permitted.  \citet{lam15}  form an initial schedule by considering each period $t$ in turn, solving a linear program with mining precedences, but ignoring processing and mining capacities, to determine which blocks to mine. A repair heuristic repeatedly selects a block, of those mined in $t$, to remove from the schedule, while ensuring that mining precedences are not violated. A variable neighbourhood descent  heuristic improves the quality of this initial schedule by repeatedly swapping the mining of two blocks in consecutive periods, shifting the mining of a block (and its successors) from one period to the next, or moving the mining of a block (and its predecessors) forward by one period. 
 
Under geological uncertainty, 
\citet{LAMGHARI2016273} compare several heuristics: the Tabu search approach of \citet{lam12}; variable neighbourhood descent \citep{lam15}; large neighbourhood search based on network flow techniques; and diversified local search. Network flow and diversified local search were found to be the most  efficient and robust. The network flow method employs two neighbourhood structures -- forward and backward -- each based on selecting blocks mined in one period to be moved to the next (or the previous) period, continually shifting the mining of blocks forward or backward in time until we achieve a new and feasible solution.  Diversified local search alternates the application of variable neighbourhood descent and the network flow strategy. 

\citet{senecal2020} apply a parallel  multi-neighbourhood form of Tabu search. The neighbourhood structures considered involve  moving the mining of a  block to a different time period, or changing its destination. Single or simultaneous applications of these moves form neighbourhoods around a current solution. The parallel algorithm maintains a pool of moves, grouped according to neighbourhood type.  Multiple threads operate in parallel,  each taking a move group from the pool, computing the result of applying each combination of moves, and maintaining a record of those resulting in the most improvement. Once the pool is empty, the best move found across threads is applied, and Tabu structures are updated. This process is repeated until a stopping criterion is met. 

 \citet{lamghari2022adaptive} define an adaptive LNS method for long-term mine planning with uncertainty. To generate an initial solution, the planning horizon is decomposed into smaller sub-problems, each solved in two stages. The first focuses on extraction decisions, and the second on material destinations. Their approach  employs multiple ways of destroying a schedule -- selecting blocks to remove -- and then repairing it -- adding blocks to the schedule. Heuristics are used to identify blocks for removal,  including: random selection; selection to reduce mining surpluses; selection based on time period mined or destination; and selection based on geological proximity and dependencies. Given a set of removed blocks, repair heuristics are then used to re-insert each block. 
Their final repair strategy -- MIP repair -- is most similar to our approach, in that it uses a MIP solver to find new values for the variables associated with removed blocks. Their approach is \textit{adaptive} as the choice of destroy/repair methods to employ at any given point in the algorithm is determined probabilistically and based on  past performance.

We consider long-term planning only in the context of open-pit mining systems. \citet{doi:10.1287/inte.2021.1087} provide a review of operations research techniques applied to underground mine planning problems. For a selection of case studies of such work, see \citet{doi:10.1287/inte.31.4.50.9669} and \citet{doi:10.1287/inte.1030.0059}.

\section{Model Formulation}\label{sec:ModelDescriptive}
Our long-term mine planning problem is defined in terms of a set of mines connected by rail to blending junctions. Each mine contains a set of pits, and each pit a set of blocks. A block contains  parcels of different material types, including high grade, low grade and waste. The optimization problem aims to find a yearly block extraction sequence that maximizes the net present value of multiple products.  Extracted material flows through a network composed of different elements -- such as crushers, stockpiles, and processing plants -- transforming ore into multiple products.

\begin{figure}[t]
\centering
\FIGURE{
    \includegraphics[width=0.8\textwidth]{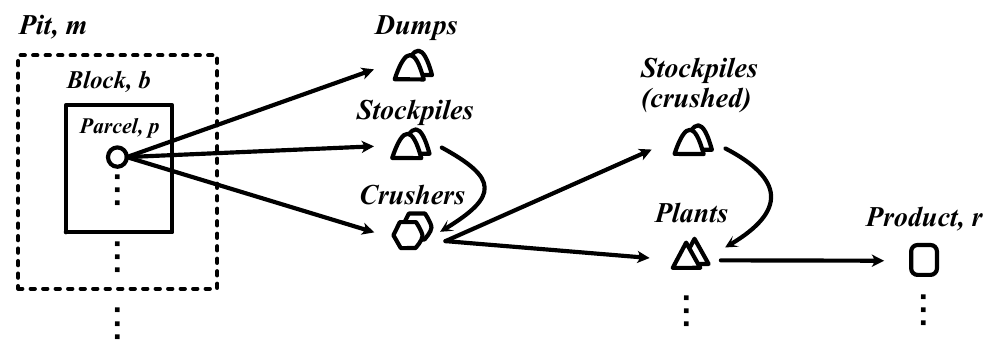}
}
{Graphical depiction of material flows for a pit in our  model.  
\label{fig:Model}}
{}
\end{figure}

This section provides a plain language description of the mathematical model capturing our long-term open-pit mine planning problem. We view a set of mines as a collection of pits, with constraints formed over sets of pits.  A flow network is defined for each pit, containing a source node denoting the pit from which blocks are extracted, and a number of destinations. These destinations capture stockpiles, crushers, processing plants, dumps, and a node for each product. Figure \ref{fig:Model} shows how material flows from sources to destinations, and between destinations. Each pit contains a set of blocks, and each block  a set of parcels representing different material types. Continuous flow variables define how much of each parcel  flows along each of the solid arrows connecting parcels to destinations (dumps, stockpiles and crushers), and between destinations (stockpiles and crushers, crushers and plants, and plants to products). The full mathematical model is provided in the appendix (\textit{Model Formulation}). Table \ref{tab:Decisions} defines the decision variables of our model.

\begin{table}[tbp]
\centering
\TABLE{Decisions present in our long-term open-pit mine planning model.\label{tab:Decisions}}
{\begin{tabular}{lp{340pt}}
\toprule
\textbf{Decision} & \textbf{Description} \\
\toprule
\textbf{Block flow} & Continuous flow variables are defined for each available pathway (source to destination nodes, and destination to product nodes) that material from each parcel may travel in each time period. 
\\[5pt]
\textbf{Block progress} & Continuous variables defining the fraction of each block extracted in each time period. \\[5pt]
\textbf{Block start} & Binary variables indicating whether extraction of a block has commenced by the end of a time period.  \\[5pt]
\textbf{Block depletion} & Binary variables indicating whether mining on a block has finished by the end of a time period. \\[5pt]
\textbf{Stockpile} & Continuous variables defining the fraction of a parcel present in a stockpile in a given time period. \\[5pt]
\textbf{Pit trigger} & Binary variables indicating whether a pit has been \textit{opened} in a given time period. Additional binary variables indicate whether the pit \textit{has been opened} in, or prior to, a given time period.\\
\bottomrule
\end{tabular}}
{}
\end{table}

\subsection{Objective}\label{sec:ObjPlain}

Our model maximizes the net present value (NPV) of a mining complex over a long-term planning horizon, capturing the total value of the investment, defined by the sum of discounted cash flows for each time period. These cash flows represent the revenue generated in each period from the sale of products minus incurred costs. In a given time period, these costs consist of capex investment costs associated with  pits opened in the  period, and operational costs. Each parcel is associated with a cost, incurred per ton extracted.

\subsection{Constraints}\label{sec:ConsPlain}

Our model contains the following key constraints.  \textit{Precedence constraints} ensure feasibility of the block extraction sequence according to the geological model and slope constraints. \textit{Blending constraints} restrict the levels of various mineral elements in each product. \textit{Capacity constraints} limit the tons of material flowing through crushers, plants, and stockpiles, and the tons of material mined in each pit. \textit{Minimum production constraints}  enforce lower bounds on the material extracted from  collections of pits.  \textit{Mass balance constraints} keep track of the fraction of each parcel present in each stockpile over time.

\section{Recovery of a Feasible Solution}\label{sec:Recovery}

Sliding windows is a popular technique for tackling large-scale time-indexed optimization problems. It   splits the planning horizon into  \textit{windows}, and solves the problem for each window in turn. Our instantiation is based on the work of \citet{cul11}, \citet{doi:10.1287/inte.2014.0737}, and \citet{o2015optimization}. We use the heuristic to find initial feasible solutions to our planning models, which are then fed into LNS. 

\begin{figure}[t]
\centering
\FIGURE
{\fbox{
    \includegraphics[width=0.8\textwidth]{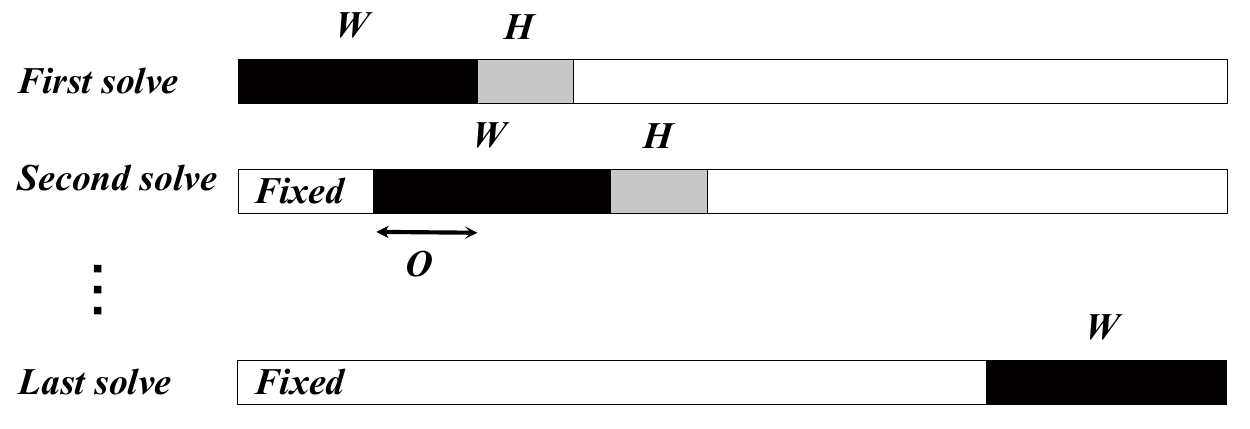}
}}
{Graphical depiction the sliding windows heuristic used to find initial feasible solutions to our long-term planning models. A solution is incrementally built by solving the problem over a window of periods $W$ and a relaxation of $H$ subsequent periods. The activities in $W - O$ of these periods, where $O$ is an overlap, are fixed, and the window moves forward. 
\label{fig:SlidingWindows}}
{}
\end{figure}

Figure \ref{fig:SlidingWindows} graphically depicts the sliding windows method. Let $W$ denote the number of time periods in each window (its width). The heuristic starts by solving the model over periods 1 to $W$, and then fixing all decisions in periods 1 to $W - O$, where $O$ is an \textit{overlap}, to their values in the resulting solution. We then move the start of our window forward to period $W + 1 - O$, and repeat the process until the last time period in our horizon has been scheduled. To reduce the likelihood of infeasibility as a result of the short-term view of each window, windows are extended with $H$ \textit{relaxed} time periods.  We then solve subproblems containing $W + H$ time periods, with all integer variables in the last $H$ of these periods becoming continuous. Note that only the activities of the first $W - O$ periods of this time frame are fixed before the window is moved forward.

\section{Large Neighbourhood Search}\label{sec:LNS}
\begin{figure}[t]
\centering
\FIGURE
{
    \includegraphics[width=0.9\textwidth]{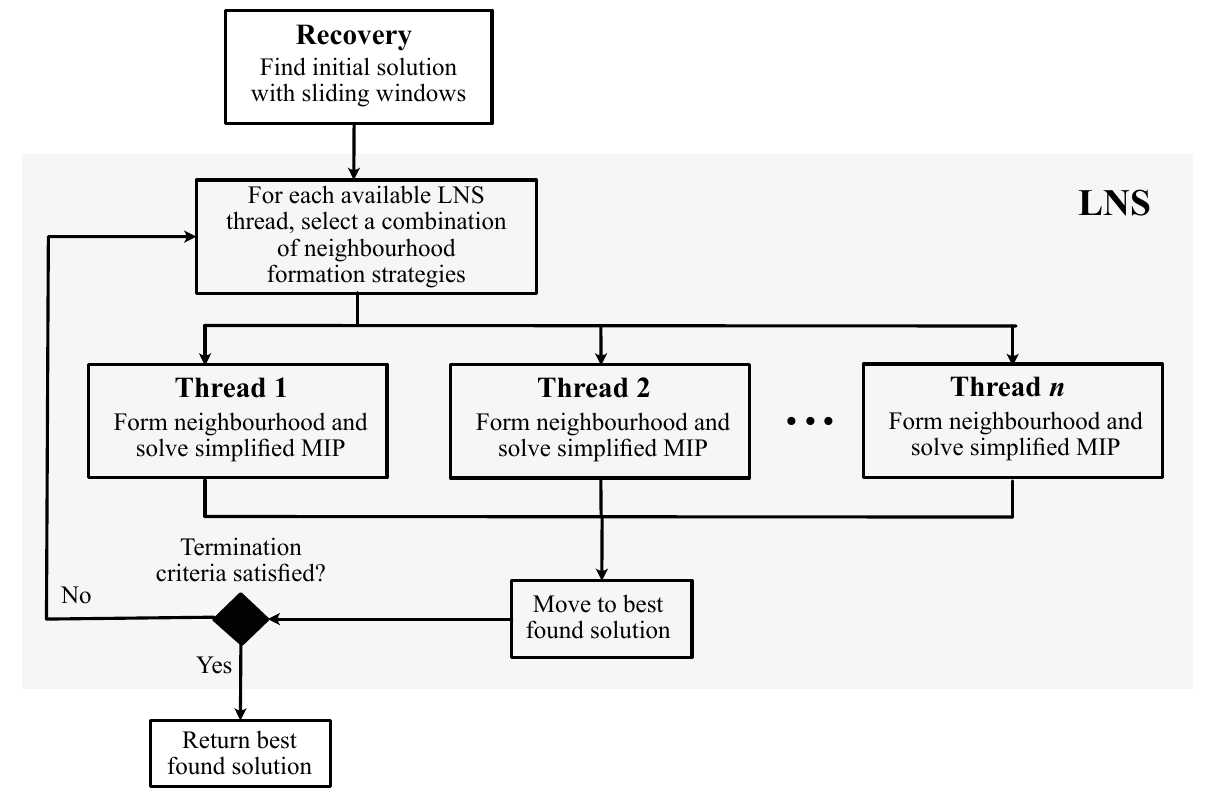}
}
{Overview of the Large Neighbourhood Search (LNS) based solution approach. Starting with an initial feasible solution, each neighbourhood formation strategy, or combination of strategies, identifies a set of variables to leave \textit{unfixed}, while fixing others, to form a simplified subproblem. \label{fig:LNSOverview}}
{}
\end{figure}

Given an initial solution formed by sliding windows, Figure \ref{fig:LNSOverview} shows how LNS is used to iteratively find improved solutions. We use a path-based neighbourhood formation process and a series of strategies to guide neighbourhood formation in each iteration. These strategies exploit characteristics of the model, and of the current best found solution. Each neighbourhood defines which variables to leave fixed to their current value, and which to leave free to take on new values, in each subproblem. We use four neighbourhood formation strategies, as outlined in Table \ref{tab:Strategies}, and described in detail in subsequent sections.

\begin{table}[t]
\centering
\TABLE{Four neighbourhood formation strategies used during LNS.
\label{tab:Strategies}}
{\begin{tabular}{lp{350pt}}
\toprule
\textsc{Strategy} &  \textsc{Description} \\
\midrule
Blending  &  Seeks to support schedule changes while  respecting blending constraints. \\[5pt]
Timing &  Supports changes to the scheduling of multiple blocks  across multiple time periods. \\[5pt]
Pit links &  Exploits the presence of minimum  production constraints. \\[5pt]
Trigger &  Seeks to improve net present value -- our objective -- by allowing key capex decisions, such as the opening or closing  of a pit, to be changed.\\
\bottomrule
\end{tabular}}
{}
\end{table}

\subsection{Path-Based Neighbourhoods}\label{sec:PBN}

Given the highly constrained nature of our model, we observe that the likelihood of finding  new, and better, solutions is low unless our neighbourhood allows changes in the extraction sequence across multiple time periods, involving blocks that span precedence chains, and that have the potential to increase NPV if the timing of their extraction changes. 
We may have a block that, if mined earlier, would increase NPV. To make this change, and retain a feasible schedule, we may need to adjust the scheduling of a precedence chain around that block. 

\begin{figure}[t]
\centering
\FIGURE
{
    \includegraphics[width=\textwidth]{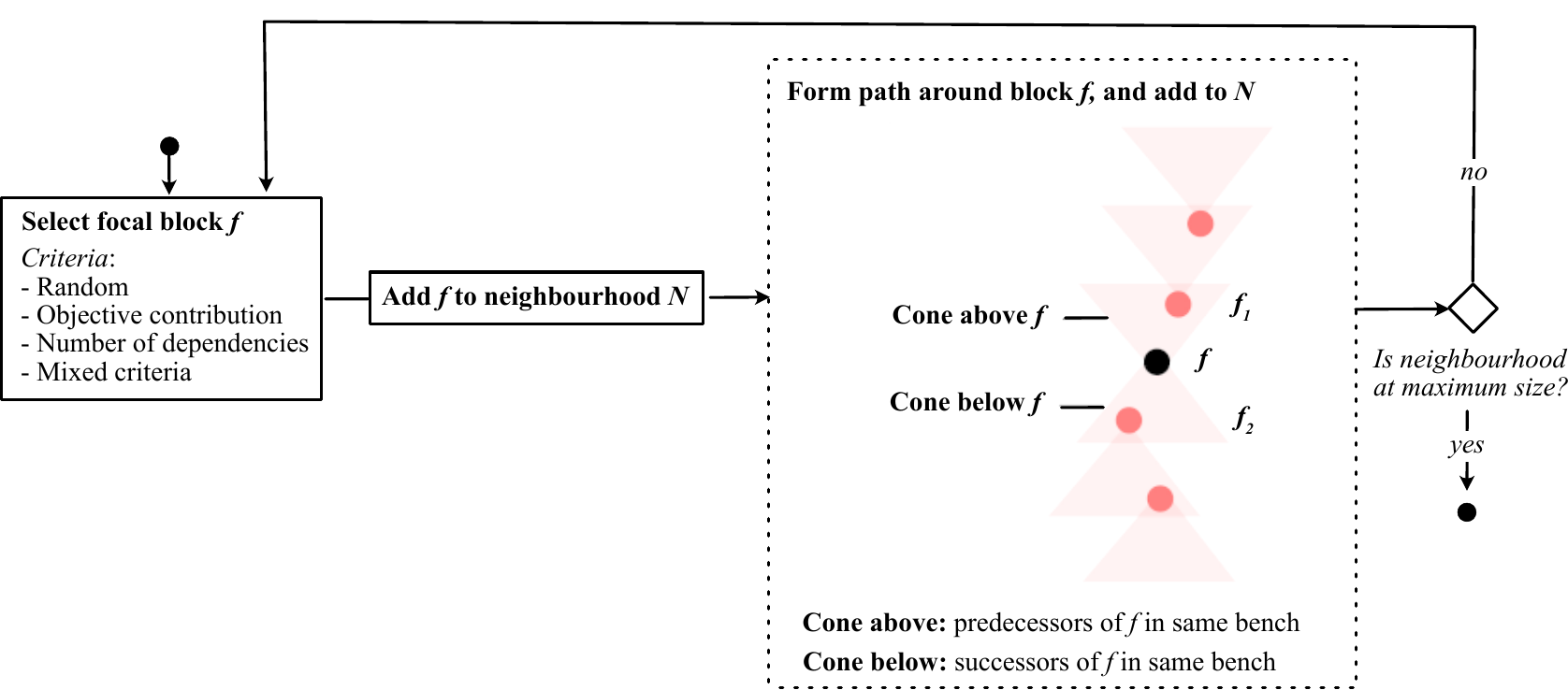}
}
{Neighbourhood formation. A selected focal block $f$ and the blocks in its restricted cone above (predecessors in the same bench) and below (successors in the same bench) form a path. The chosen block selection criterion is used to select two new focal blocks, one from the restricted cone above $f$, $f_1$, and below, $f_2$. Blocks in the restricted cone above $f_1$, and below $f_2$, are added to the path. This process is repeated, selecting a new focal block in the cone above $f_1$, and in the cone below $f_2$, and continues until the cones above and below selected focal blocks are empty.   \label{fig:LNSFormation}}
{}
\end{figure}

Figure \ref{fig:LNSFormation} shows how we form a neighbourhood of blocks. We form neighbourhoods by repeatedly selecting a \textit{focal} block from the set of all blocks, forming a \textit{path} around that block, and adding the blocks in that path to our neighbourhood. This process is repeated until our neighbourhood reaches a maximum size.

 To form neighbourhoods with blocks that are likely to have an impact on NPV, we select focal blocks randomly on the basis of a weighted distribution. Each block is assigned a weight, designed to represent its value for the purposes of scheduling. We consider several methods for weighting blocks: random; objective; minimum dependencies; and mixed criteria. These methods are described in Table \ref{tab:BlockSelection}. 

\begin{table}[t]
\centering
\TABLE{Block selection strategies used in neighbourhood formation.
\label{tab:BlockSelection}}
{\begin{tabular}{p{100pt}p{300pt}}
\toprule
\textsc{Strategy} & \textsc{Description} \\
\midrule
\textbf{Random (RAND)} & Focal blocks are selected randomly from the set of all blocks.\\
\hline
\textbf{Objective (OBJ)} & The contribution of a block to the objective function is used as its weight. We use the iron (Fe) content of the block as a proxy for this contribution.\\
\hline
\textbf{Minimum Dependencies (MD)} & Blocks with fewer dependencies are given higher weights than those with more. The weight of block with $n$ predecessors and successors, is set to $1/n$ for $n > 0$ and 0 otherwise. The idea  is to bias selection toward potential bottleneck blocks that, when scheduled, have a substantial influence on the overall schedule.  \\
\hline
\textbf{Mixed Criteria (MIX)} & Each LNS solving thread is assigned one of the three weighting methods to use -- random is assigned to the first thread, objective to the second, and minimum dependencies to the third, with the pattern repeated across the remaining threads. \\
\bottomrule
\end{tabular}}
{}
\end{table}
\ignore{
 \textbf{Path Formation} To form a path around a block $f$, we first consider the blocks in its \textit{restricted cone above}, $CA(f)$. These are the predecessors of $f$ that belong to the same bench as $f$. Both $f$, and $CA(f)$, are added to $p(f)$. We then consider the blocks in its \textit{restricted cone below}, $CB(f)$, which are also added to $p(f)$. To expand this path, we randomly select a block $b$ in $CA(f)$ according to their values (see \textit{Block Selection}). Block $b$ and $CA(b)$ are added to $p(f)$. At the same time, we select a block $b'$ in $CB(f)$, adding both $b'$ and $CB(b')$ to $p(f)$. We continue this process of expanding our path upwards, and downwards, until the cones above and below our selected blocks are both empty, or we have reached the maximum neighbourhood size, $\bar{N}$. 
Figure \ref{fig:PathExample} visually depicts the formation of a path-based neighbourhood.

 \begin{figure}[tbp]
\centering
\FIGURE
{\includegraphics[width=5cm]{Figures/PathBasedN.png}}
{The formation of a \textit{path} of blocks, around a focal block $f \in \mathcal{B}$.\label{fig:PathExample}}
{}
\end{figure}}

For a neighbourhood, $N$, we form a restricted MIP as follows. For each block that \textit{does not belong to} $N$, we fix the mining decisions relating to that block to their values in the current best solution. The decisions we consider fixing for a block are its start and depletion binaries, and continuous flow variables. 

\subsection{Neighbourhood Formation Strategies}\label{sec:NStrategies}
This section defines several variations on the previously described method of forming path-based neighbourhoods, each designed to exploit different aspects of the model structure, and current best solution, to maximize the likelihood of finding a better solution. While these strategies exploit specific features of our model, these features are typical of real-world long-term mine planning models. Where multiple strategies are used, each LNS thread is set up to apply varying combinations of these strategies.

\textbf{Blending}  Each block makes a certain contribution to a range of blending constraints, for each of the available paths along which material from the block can travel in each time period. This contribution may be positive, negative, or zero. These constraints are satisfied if the sum of these contributions is greater than zero. Using the Blending strategy, we keep track of the total contribution of blocks already in the neighbourhood being formed to each of the blending constraints. While the blocks in the neighbourhood will be mined across different time periods, we consider an approximation in which we imagine the total contribution of the neighbourhood to each blending constraint, as if it were being mined all at once.

The Blending strategy restricts focal block selection to the set of blocks with  non-zero weights, and that do not make a negative contribution to any blending constraint where the total contribution of the neighbourhood is negative. This is  achieved by setting the weight of blocks that do not satisfy this criterion to zero. To form a path around a focal block, the process described in the \textit{Path-Based Neighbourhoods} section is followed, with no restrictions on which blocks are selected to extend the path.  

\textbf{Timing} When forming a neighbourhood, once we select a focal block $f$ that is completely mined within the planning horizon, we restrict all subsequent focal block selections to the set of blocks with non-zero weights and that either start being extracted, or are mined to completion, while $f$ is being mined. The idea  is to form a neighbourhood with subsets of blocks mined at similar times.

\textbf{Pit Links}  To find  improved solutions, we need to be able to  make changes to an existing schedule while still satisfying minimum  production constraints over sets of pits. 
Under the Pit Links strategy, blocks are initially weighted according to one of the methods in Table \ref{tab:BlockSelection}.  After adding the first path of blocks  to the neighbourhood, these values are adjusted for subsequent focal block selections.
We set the weight of a block to zero if it does not belong to a pit that is either \textit{covered} by the current neighbourhood, or linked, by a minimum production constraint, to a pit covered by the neighbourhood. A neighbourhood \textit{covers} a pit if it contains a block that belongs to that pit.  
If the weights of all remaining blocks are zero, and we have not reached the maximum neighbourhood size, we reset all block weights to their original values.

\textbf{Trigger} The choice of when, or if, to open a pit, incurring capex costs, is likely to have a significant impact on NPV. In our model, binary variables, denoted \textit{pit triggers}, are associated with the opening of each pit. We can associate each pit trigger variable  with a subset of blocks. The scheduling of these blocks  determines the timing of these  decisions. 
Under the Trigger strategy, we do not form path-based neighbourhoods. Instead, we randomly select a trigger variable, and include its associated blocks in the neighbourhood being constructed. We repeat this process until we reach a desired neighbourhood size.

\subsection{Solving Strategies}\label{sec:SolvingStrategies}
We employ a number of additional strategies to improve the performance of the overall algorithm. 

\textbf{Earliest Starts}
The computation of earliest extraction periods for each block in a mine planning model is a well known strategy for reducing model complexity and solve times. For a given block, $b$, we compute the total tons of material that would have to be extracted, on the basis of mining precedence constraints, before extraction of $b$ can begin. For a given time period, $t$, if this value exceeds the total cumulative mining capacity across periods 1 to $t$, then we can fix all mining related decision variables for $b$ -- block start, depletion, and flow variables -- in period $t$,  to zero. 

\textbf{Relaxation Induced Neighbourhood Search (RINS)} RINS is a neighbourhood construction approach that utilises a current incumbent solution and information obtained by solving the linear relaxation of a MIP \citep{danna2005exploring}. The idea is that variables that take on the same value in the incumbent, and in the MIPs linear relaxation, will likely take on the same value in a good solution. We borrow this idea to allow larger neighbourhood sizes, in terms of numbers of blocks, while retaining reasonable subproblem solve times. 
For a given LNS subproblem, we first solve its linear relaxation. We fix each integer variable that takes on the same value in the current best solution (the incumbent), and in the solution to the linear relaxation, to that value. The subproblem is then re-solved.

\textbf{Fixing variables}
We consider two different levels of fixing: S/D; and S/D + F. Under S/D, we fix only the block start and depletion variables for blocks that do not belong the neighbourhood $N$. Under S/D + F, we \textit{additionally} fix the flow variables for blocks that \textit{do not belong to the same pit as any block in} $N$. 
We anticipate that for most blocks, we will not be able to radically change the time at which they are mined in a single iteration of LNS. We define an \textit{unfix window}, $UW$, to control the size of each subproblem. 
For any block  whose decision variables are to be unfixed, they are only unfixed in a time range starting $UW$ periods before the block is first mined, to $UW$ periods after the block is depleted (or the end of the horizon if it is never depleted). If the block is never mined, its decision variables are left unfixed across the entire horizon.

\textbf{Termination Criteria}
We define two criteria to determine when to terminate LNS:  a \textit{time limit} on run time;  and failure to maintain a desired \textit{improvement rate}. In the latter, LNS terminates after a minimum number of iterations have been performed, and the rate of solution improvement -- the sum of percentage improvements across all iterations thus far, divided by the number of iterations -- falls below a threshold. 

\section{Results}\label{sec:Evaluation}
All experiments have been conducted on a machine with an Intel Xeon Platinum 8176 chip (2.1GHz), and 1TB of RAM. All MIPs are solved with Gurobi 9.11 using the following parameter settings: NumericFocus (1); MIPFocus (1); PrePasses (2); AggFill (200); and Heuristics (0.25). All other parameters are set to their default values. Through informal experimentation, this parameter setting was found to elicit the best performance from Gurobi when solving the full MIP models of our two problem instances, and performed well when solving subproblem MIPs during LNS. In all reported runs of LNS, neighbourhoods are explored over 10 parallel threads. Earliest starts for each block are used when solving full MIP models and all LNS subproblems, across all experiments.  The sliding windows heuristic and LNS were implemented in C++. 

An evaluation of the sliding windows heuristic across varying  window size ($W \in \{1\ldots7\}$), relaxed horizon length ($H \in \{0\ldots4\}$), and overlap ($O \in \{0\ldots2\}$), can be found in the appendix (\textit{Detailed Results}). We report the quality of solutions found, in terms of their gap to either a known optimal solution or a best known upper bound, for the Oyu Tolgoi and Pilbara models. Increasing the overlap, or window width, generally results in a better initial solution at the cost of increased run times. Increasing the horizon length also has the capacity to increase solution quality, albeit less reliably. Controlling $W$  appears to be the most effective method for achieving a desired trade off between solution quality and solve time.

The parameters $W = 5$, $H = 2$, and $O = 1$ result in a reasonable initial starting solution for Oyu Tolgoi, with an optimality gap of 27.1\%, found in 98s. Across the parameter combinations considered for the Pilbara, $W = 4$, $H = 0$, and $O = 1$ finds an initial solution with the smallest gap, 8.2\%, in 82 minutes.

For both case studies, we start with a baseline parameter configuration. Through informal experimentation, these settings were found to be the most performant of those considered. We then examine the impact of adjusting each parameter in turn to understand the relative value of different options.

\subsection{Oyu Tolgoi}

We apply LNS to an initial solution with an optimality gap of 27.1\%, found with sliding windows $W = 5$, $H = 2$, and $O = 1$. Table \ref{tab:OTNSizeVary}, in the appendix (\textit{Detailed Results}), reports the average quality of solutions found by LNS  across varying neighbourhood sizes (50 to 150), focal block selection methods (OBJ, MD, RAND, and MIX), variable fixing strategies (S/D and S/D + F), use of RINS as a solving strategy (yes/no), and neighbourhood formation strategies (Timing and Blending, Timing only, Blending only, and None). We report solution quality in terms of the gap to the best known upper bound on the objective, found after solving the full model with Gurobi for 148 hours, warm started with the same initial solution used as the starting point for LNS.  We set the unfix window $UW$ to the length of the planning horizon. Oyu Tolgoi has blending constraints, but no minimum production constraints. Consequently, both the Timing and Blending strategies are applicable. 
 Across all LNS runs on the Oyu Tolgoi model, MIPs are terminated at a gap of 0.05\%, or after 600s, and LNS is terminated after 12 hrs or when the improvement rate falls below 1\%.

We consider a baseline setting with a neighbourhood of 100 blocks, objective (OBJ) based focal block selection, S/D variable fixing (block start and depletion variables), and the use of RINS and the Timing and Blending neighbourhood formation strategies (Table \ref{tab:OTNSizeVary}, entry 2). This setting produces solutions with an average gap of 10.3\% in 258 minutes, across 10 runs. 
Figure \ref{fig:OTLnsVsDefault} compares the performance of LNS in this setting against solving a warm-started Oyu Tolgoi model directly with Gurobi. 
While the full Gurobi solve eventually results in a superior solution, it requires substantially more time.

\begin{figure}[!t]
    \centering
    \FIGURE
    {
\begin{tikzpicture}
\begin{axis}[restrict x to domain=0:1500, restrict y to domain=0:60,
	legend style = {draw=none,at={(0.95,0.95)},anchor=north east},
    title={Oyu Tolgoi -- LNS vs Gurobi},
    xmode=log,
    log ticks with fixed point,
	scaled ticks=false,
    xtick pos=left,
    ytick pos=left,
	xlabel = {Time (mins)}, 
    xtick ={10, 100,   500,  1500},
	ylabel = {Gap to Optimal (\%)},
	height=7cm, width=10cm]
	\addplot+[solid,black,mark=none] table [x = TIME, y = GAP]{PaperData/OT_FULL.csv};
	\addlegendentry{Gurobi}
	\addplot+[solid,blue,mark=none] table [x = TIME, y = GAP]{PaperData/ot_mins_s1.dat};
	\addlegendentry{LNS}
	\addplot+[solid,blue,mark=none] table [x = TIME, y = GAP]{PaperData/ot_mins_s2.dat};
	\addplot+[solid,blue,mark=none] table [x = TIME, y = GAP]{PaperData/ot_mins_s3.dat};
	\addplot+[solid,blue,mark=none] table [x = TIME, y = GAP]{PaperData/ot_mins_s4.dat};
 	\addplot+[solid,blue,mark=none] table [x = TIME, y = GAP]{PaperData/ot_mins_s5.dat};
   	\addplot+[solid,blue,mark=none] table [x = TIME, y = GAP]{PaperData/ot_mins_s6.dat};
     	\addplot+[solid,blue,mark=none] table [x = TIME, y = GAP]{PaperData/ot_mins_s7.dat};
       	\addplot+[solid,blue,mark=none] table [x = TIME, y = GAP]{PaperData/ot_mins_s8.dat};
         	\addplot+[solid,blue,mark=none] table [x = TIME, y = GAP]{PaperData/ot_mins_s9.dat};
 	\addplot+[solid,blue,mark=none] table [x = TIME, y = GAP]{PaperData/ot_mins_s10.dat}; 
\end{axis}
\end{tikzpicture}}
    {\textbf{Oyu Tolgoi}: Optimality gap over time across 10 LNS runs,  using the baseline setting and 10 LNS threads. Run time (logarithmic scale) is inclusive of  sliding windows, forming an initial feasible solution with a 27.1\% gap. LNS is compared against a full solve with Gurobi, warm started with the same starting solution. All ten LNS solves reach a gap of under 11\% within 5.5 hours of solving.    \label{fig:OTLnsVsDefault}}
    {}
\end{figure}

Table \ref{tab:OTBlockSize} (left) shows how reducing  \textit{neighbourhood size} from 100 (baseline) to 50 blocks substantially reduces both average LNS run times and solution quality. When MIP subproblem sizes become too large to be effectively solved in the 600s time limit (150 blocks), solution quality also degrades. Table \ref{tab:OTBlockSize} (right) shows that varying the \textit{focal block selection method} does not result in substantial differences in solution quality, while  objective-based selection reduces average solve times by 18 to 24\%. 

\begin{table}[t]
\centering
\TABLE{\textbf{Oyu Tolgoi}: Impact of neighbourhood size (left) and focal block selection method (right) on the average time (min) spent by LNS, and average quality of solutions (\% gap), across 10 runs. 
\label{tab:OTBlockSize}}
{\begin{tabular}{c||c|c}
\toprule
     &   Avg Time    & Avg Gap \\
  $N$   & (Min, Max) (min) & (Min, Max)  (\%) \\
 \midrule
50 &  38.3 (25,50)   & 17.5 (12.5,26.0)    \\
 100 &   258.3, (161.7,315)   &    10.3 (9.1,11.3) \\
150 &  378.3, (343.3,461.7)   &   10.5 (9.0,11.5) \\
 \bottomrule
 \multicolumn{2}{c}{}\\
\end{tabular}\label{tab:OTBlockSize} 
\quad
\begin{tabular}{c||c|c}
\toprule
  Focal   &   Avg Time    & Avg Gap \\
  block   & (Min, Max) (min) & (Min, Max)  (\%) \\
 \midrule
 OBJ &  258.3, (161.7,315)   &    10.3 (9.1,11.3) \\
 MD &  323.3, (290,353.3)   &    9.9 (7.8,11.2)\\
 RAND &  338.3, (300,391.7)   &    9.6 (7.9,10.8)\\
 MIX &  315, (280,353.3)   &   10.0 (8.4,10.7)\\
 \bottomrule
\end{tabular}\label{tab:OTFocalBlock}
}
{}
\end{table}

 The additional \textit{fixing of flow variables} (S/D + F) results in a similar average run time and solution quality, 252 minutes and a gap of 10.7\%, to the baseline (Table \ref{tab:OTNSizeVary}, entry 3).  Without the \textit{use of RINS}, the average optimality gap of solutions increases to 13\%,  and average run times reduce to 205 minutes (Table \ref{tab:OTNSizeVary}, entry 4).  Without RINS, MIP solves are more likely to be cut off at the time limit of 600s before reaching a desired gap, leading to smaller improvement rates over time, and earlier termination.

\begin{table}[h]
\centering
\TABLE{\textbf{Oyu Tolgoi}: Impact of neighbourhood formation methods on the average time (min) spent by LNS, and average quality of solutions (\% gap), across 10 runs. 
\label{tab:OTStrategies}}
{\begin{tabular}{c||c|c}
\toprule
     &   Avg Time    & Avg Gap \\
  Strategies   & (Min, Max) (min) & (Min, Max)  (\%) \\
 \midrule
\textit{Timing} \& \textit{Blending} &  258.3, (161.7,315)   &   10.3 (9.1,11.3) \\
 \textit{Timing} only &   133.3, (83.3,291.7)   &   11.4 (10.2,13.1) \\
 \textit{Blending} only &  285, (211.7,355)   &   10.5 (8.8,12.0)  \\
 None &  93.3, (28.3,176.7)   &   13.9 (11.1,26.0)\\
 \bottomrule
\end{tabular}}
{}
\end{table}

Table \ref{tab:OTStrategies} shows the impact of varying the combination of \textit{neighbourhood formation strategies} used on average solution quality and LNS run times. The Blending strategy is the most impactful in this context, in terms of reducing the average optimality gap of solutions. While using either of the applicable strategies (Timing and Blending) alone is better than using none, the best performance is achieved when using both.

\subsection{Pilbara}

We apply LNS to an initial solution with an optimality gap of 8.2\%, found with sliding windows in 82 minutes using $W = 4$, $H = 0$, and $O = 1$. Table \ref{tab:PilbaraLNS}, in the appendix (\textit{Detailed Results}), reports the average quality of solutions found by LNS across varying focal block selection methods (OBJ, MD, RAND, and MIX), variable fixing methods (S/D and S/D + F), use of RINS (yes/no), choice of unfix window ($UW \in [2,4]$), and neighbourhood formation strategies used (None, All, Trigger only, Timing and Blending, Timing and Pit Links, All except Blending, and All except Pit Links). We use a fixed neighbourhood size of 200 blocks across all experiments.  The Pilbara model contains blending constraints, capex decisions, and minimum  production constraints. All neighbourhood formation strategies in Table \ref{tab:Strategies} are applicable.

For each parameter setting, Table \ref{tab:PilbaraLNS} reports the average quality of solutions found by LNS, in terms of their gap to the best known upper bound on the optimal objective, over 10 differently seeded runs, each utilising 10 threads of computation. LNS MIPs were solved to a gap of 0.05\%, or a time limit of 400s. Each LNS run was terminated after 12 hours.
When warm-started with the initial feasible solution used by LNS, Gurobi was unable to improve upon the solution after 148 hours of solve time. We compare the best upper bound on the optimal objective found at this point with the quality of solutions found by LNS.

\begin{figure}[!t]
    \centering
    \FIGURE
    {\begin{tikzpicture}
\begin{axis}[restrict x to domain=0:800, restrict y to domain=0:12,
	legend style = {draw=none,at={(0.9,0.85)},anchor=north east},
    title={Pilbara -- LNS, 10 Runs},
    title style={yshift=-0.4ex,},
	scaled ticks=false,
    xtick pos=left,
    ytick pos=left,
	xlabel = {Time (min)}, 
    xtick ={0, 200,  400,  600,  800},
	ylabel = {Gap to Optimal (\%)},
	height=7cm, width=10cm]
 	\addplot+[solid,black,mark=none] coordinates {(0,9.32) (800,9.32)};
	\addlegendentry{Gurobi}
	\addplot+[solid,blue,mark=none] table [x = TIME, y = GAP]{PaperData/pilbara_best_lns_mins_s1.dat};
	\addlegendentry{LNS}
	\addplot+[solid,blue,mark=none] table [x = TIME, y = GAP]{PaperData/pilbara_best_lns_mins_s2.dat};
	\addplot+[solid,blue,mark=none] table [x = TIME, y = GAP]{PaperData/pilbara_best_lns_mins_s3.dat};
	\addplot+[solid,blue,mark=none] table [x = TIME, y = GAP]{PaperData/pilbara_best_lns_mins_s4.dat};
 	\addplot+[solid,blue,mark=none] table [x = TIME, y = GAP]{PaperData/pilbara_best_lns_mins_s5.dat};
 	\addplot+[solid,blue,mark=none] table [x = TIME, y = GAP]{PaperData/pilbara_best_lns_mins_s6.dat};
 	\addplot+[solid,blue,mark=none] table [x = TIME, y = GAP]{PaperData/pilbara_best_lns_mins_s7.dat};
 	\addplot+[solid,blue,mark=none] table [x = TIME, y = GAP]{PaperData/pilbara_best_lns_mins_s8.dat};
 	\addplot+[solid,blue,mark=none] table [x = TIME, y = GAP]{PaperData/pilbara_best_lns_mins_s9.dat};
 	\addplot+[solid,blue,mark=none] table [x = TIME, y = GAP]{PaperData/pilbara_best_lns_mins_s10.dat};
\end{axis}
\end{tikzpicture}}
{\textbf{Pilbara}: Optimality gap over time across 10 LNS runs, using the baseline setting and 10 LNS threads. Run time is inclusive of  sliding windows, forming an initial solution with a 8.2\% gap. LNS is compared against a full solve with Gurobi, warm started with the same starting solution. \label{fig:PilbaraLNS_plots}}
{}
\end{figure}

We consider a baseline setting in which all of the neighbourhood formation strategies have been used, objective-based focal block selection, flow variables fixed (S/D + F), no RINS applied during solving, and an unfix window of $UW = 2$ (Table \ref{tab:PilbaraLNS}, entry 8). This setting  produces solutions with an average gap of 4.5\% over  12 hours. Figure \ref{fig:PilbaraLNS_plots} compares the performance of LNS, using the baseline setting and 10 LNS threads, against solving a warm-started Pilbara model directly with Gurobi.

Our results suggest that the Trigger \textit{neighbourhood formation strategy}, which considers key capex decisions over time, is most effective for the Pilbara model. Modifying the baseline setting to only use the Trigger strategy, in place of the entire set, does not change the average quality of the solutions found (Table \ref{tab:PilbaraLNS}, entry 5). The difference in solution quality when using the Trigger strategy versus none at all (average gaps of 4.5 versus 4.7) was small but statistically significant (Wilcoxon signed rank test applied to the two populations of solutions: $p = 0.027$). Using only Timing and Blending, or only Timing and Pit Links, both increase the average optimality gap of solutions to 4.7 and 4.6, respectively (Table \ref{tab:PilbaraLNS}, entries 6-7). 

Modifying the baseline to use an alternate \textit{focal block selection method}, in place of objective-based, does not change the average quality of the solutions found (Table \ref{tab:PilbaraLNS}, entries 11-13).
Use of the \textit{RINS strategy} was not helpful for the Pilbara. The cost of two solves for each subproblem under  RINS  (an LP relaxation, and then a restricted MIP) was higher than a single solve of the original subproblem. Using RINS increases the average gap of the solutions found to 4.8\%, over the baseline of 4.5\% (Table \ref{tab:PilbaraLNS}, entry 3). Increasing the \textit{unfix window} from 2 to 4, leaving more variables unfixed in each subproblem, also increases the average gap of solutions to 4.9\% (Table \ref{tab:PilbaraLNS}, entry 4). 
For Oyu Tolgoi, there were no substantial differences in the quality of the solutions obtained with the two \textit{variable fixing strategies} (S/D and S/D + F). For the Pilbara, we need to fix the additional flow variables to maintain manageable subproblem sizes. Modifying the baseline setting to fix only the block start and depletion variables (S/D) increases the average optimality gap of solutions to 8.2\% (Table \ref{tab:PilbaraLNS}, entry 2). In this case, we have not improved upon the initial solution.

\section{Impact of the new scheduling method}

Rio Tinto mine planners conduct long-term planning on a quarterly basis, generating long-term schedules that will subsequently be used to inform short to medium-term planning. There are many uncertainties in the input data to the planning process and various options to consider, ranging from decisions on when to build new infrastructure to which product strategies to pursue. Different product strategies will explore different targets on the production of each product. The way Rio Tinto sells ore in the market can affect prices, and data provided by multiple schedules formed under different scenarios can be used to inform price estimation. Ideally, long-term planners should construct a scenario matrix containing multiple schedules generated with different combinations of options.  

Generating schedules in a reasonable time frame, for large systems like the Pilbara, was previously only possible with heuristics like sliding windows. To form reasonable schedules, sliding windows would be run with longer window lengths and overlap parameter values, resulting in longer solve times. From the perspective of the iron ore commodity group in particular, the 7-10 day solve time to find a single schedule has been problematic. Mine planners are required to work to deadlines that, in conjunction with these long solve times, allow them to generate only a single schedule to provide to the business, exploring a single product strategy or scenario. 
The LNS technique uses sliding windows under more aggressive parameters, forming an initial schedule quickly, which is then iteratively improved by solving simpler subproblems. This approach reduces the 7-10 day solve time for the Pilbara to a single day, or less, allowing more scenarios to be investigated and, as a result, better data provided to the business for decision-making. 

The faster solve times achievable with LNS allow mine planners to increase the complexity of their models. The rail network connecting sites within the Pilbara system, for example, is considered to be a scheduling bottleneck. The configuration of the rail network has an impact on what the optimal production schedules will be across sites, and this impact is expected to increase with future expansions. Inclusion of the rail network in the integrated Pilbara model becomes feasible with faster solving methods, providing a more holistic view of the system during optimization. For many models, faster solve times allow additional logic related to truck fleets and stockpiling to be incorporated. This allows planners to optimize for fleet replacements with electric trucks, for example. For large models such as the Pilbara, the inclusion of all truck fleets and stockpiles in the model was previously intractable. Improved solutions have been possible with LNS due to the ability to include such detail.  

\section{Concluding Remarks}\label{sec:Discussion}

We have presented a Large Neighbourhood Search  based method for solving long-term open-pit mine planning problems, demonstrating the effectiveness of the approach on two real world mining projects of varying complexity.  Our  method is able to find, within hours, near optimal solutions to planning problems in substantially reduced time frames than with  an off-the-shelf solver. Reasonable solve times are important  where multiple plans must be generated over a range of scenarios in order to make long-term strategic decisions.

Maintaining modest subproblem sizes while maximising the likelihood of finding new and improved solutions was a key challenge in the development of our approach. The key insight was that the features of the model itself -- what decisions were likely to have the most impact on net present value, and what constraints were likely to be the hardest to satisfy -- needed to be considered when forming neighbourhoods. Second, it was important to find opportunities to reduce subproblem size that did not materially reduce the chances of finding new and improved solutions. Fixing the continuous flow variables for blocks in pits not covered by a neighbourhood reduced subproblem sizes while having little impact on solution quality. Similarly, leaving decisions for neighbourhood blocks unfixed for some but not all time periods was particularly beneficial for the larger Pilbara model. For LNS to be successful, we needed to form subproblems that both (i) had a high chance of admitting a different, and better, solution than we started with, and (ii) could be effectively solved within a time limit to reach these solutions.  

Our analysis suggests several guidelines that can be used when faced with a new model. 
 As model complexity increases, strategies that keep subproblem size and solve times manageable become important. The use of the RINS strategy shifts from being helpful to detrimental, while restricting subproblem size by fixing more variables (with S/D + F and a small unfix window, $UW$) becomes important. For different models, certain neighbourhood formation strategies may be more beneficial than others. However, using all applicable strategies is unlikely to be detrimental. The choice of focal block selection method is unlikely to make a material difference in performance.

The presented scheduling method has allowed Rio Tinto to substantially improve upon the solutions provided by  heuristics like sliding windows, providing value insights on a scale of millions of dollars. In future work, these methods can be developed further by incorporating adaptivity in the selection of neighbourhood formation strategies across available threads, biasing selection based on their prior performance. In addition, where there is disparity in the solve times of subproblems across threads, a degree of asynchronicity can be considered to allow new neighbourhoods to be explored on some threads, while waiting for subproblem solves to complete on others.

\bibliographystyle{apalike} 
\bibliography{references}

\begin{thebibliography}{}

\bibitem[Amaya et~al., 2009]{ama09}
Amaya, J., Espinoza, D., Goycoolea, M., Moreno, E., Prevost, T., and Rubio, E.
  (2009).
\newblock A {S}calable {A}pproach to {O}ptimal {B}lock {S}cheduling.
\newblock In {\em APCOM}, pages 567--575.

\bibitem[Askari-Nasab et~al., 2011]{ask11JOMS}
Askari-Nasab, H., Pourrahimian, Y., Ben-Awuah, E., and Kalantari, S. (2011).
\newblock Mixed {I}nteger {L}inear {P}rogramming {F}ormulations for {O}pen
  {P}it {P}roduction {S}cheduling.
\newblock {\em Journal of Mining Science}, 47:338--359.

\bibitem[Carlyle and Eaves, 2001]{doi:10.1287/inte.31.4.50.9669}
Carlyle, W.~M. and Eaves, B.~C. (2001).
\newblock Underground planning at stillwater mining company.
\newblock {\em Interfaces}, 31(4):50--60.

\bibitem[Chicoisne et~al., 2012]{chic12}
Chicoisne, R., Espinoza, D., Goycoolea, M., Moreno, E., and Rubio, E. (2012).
\newblock A {N}ew {A}lgorithm for the {O}pen-{P}it {M}ine {P}roduction
  {S}cheduling {P}roblem.
\newblock {\em Operations Research}, 60(3):517--528.

\bibitem[Chowdu et~al., 2022]{doi:10.1287/inte.2021.1087}
Chowdu, A., Nesbitt, P., Brickey, A., and Newman, A.~M. (2022).
\newblock Operations research in underground mine planning: A review.
\newblock {\em INFORMS Journal on Applied Analytics}, 52(2):109--132.

\bibitem[Cullenbine et~al., 2011]{cul11}
Cullenbine, C., Wood, R.~K., and Newman, A. (2011).
\newblock {A} sliding time window heuristic for open pit mine block sequencing.
\newblock {\em Optimization Letters}, 5:365--377.

\bibitem[Danna et~al., 2005]{danna2005exploring}
Danna, E., Rothberg, E., and Pape, C.~L. (2005).
\newblock Exploring relaxation induced neighborhoods to improve mip solutions.
\newblock {\em Mathematical Programming}, 102(1):71--90.

\bibitem[Epstein et~al., 2012]{eps12}
Epstein, R., Goic, M., Weintraub, A., Catalan, J., Santibanez, P., Urrutia, R.,
  Cancino, R., Gaete, S., Aguayo, A., and Caro, F. (2012).
\newblock Optimizing {L}ong-{T}erm {P}roduction {P}lans in {U}nderground and
  {O}pen-{P}it {C}opper {M}ines.
\newblock {\em Operations Research}, 60:4--17.

\bibitem[Fathollahzadeh et~al., 2021]{Fathollahzadeh2021-jo}
Fathollahzadeh, K., Asad, M. W.~A., Mardaneh, E., and Cigla, M. (2021).
\newblock Review of solution methodologies for open pit mine production
  scheduling problem.
\newblock {\em Int. J. Min. Reclam. Environ.}, 35(8):564--599.

\bibitem[Gershon, 1987]{gersh87}
Gershon, M. (1987).
\newblock Heuristic approaches for mine planning and production scheduling.
\newblock {\em International Journal of Mining and Geological Engineering},
  5(1):1--13.

\bibitem[Goodfellow and Dimitrakopoulos, 2016]{Goodfellow2016-uk}
Goodfellow, R.~C. and Dimitrakopoulos, R. (2016).
\newblock Global optimization of open pit mining complexes with uncertainty.
\newblock {\em Appl. Soft Comput.}, 40:292--304.

\bibitem[Kuchta et~al., 2004]{doi:10.1287/inte.1030.0059}
Kuchta, M., Newman, A., and Topal, E. (2004).
\newblock Implementing a production schedule at lkab's kiruna mine.
\newblock {\em Interfaces}, 34(2):124--134.

\bibitem[Lambert et~al., 2014]{lam14}
Lambert, W.~B., Brickey, A., Newman, A.~M., and Eurek, K. (2014).
\newblock Open-{P}it {B}lock-{S}equencing {F}ormulations: {A} {T}utorial.
\newblock {\em Interfaces}, 44(2):127--142.

\bibitem[Lamghari, 2017]{Lamghari2017-xk}
Lamghari, A. (2017).
\newblock Mine planning and oil field development: A survey and research
  potentials.
\newblock {\em Math. Geosci.}, 49(3):395--437.

\bibitem[Lamghari and Dimitrakopoulos, 2012]{lam12}
Lamghari, A. and Dimitrakopoulos, R. (2012).
\newblock A diversified tabu search approach for the open-pit mine production
  scheduling problem with metal uncertainty.
\newblock {\em European Journal of Operational Research}, 222(3):642--652.

\bibitem[Lamghari and Dimitrakopoulos, 2016]{LAMGHARI2016273}
Lamghari, A. and Dimitrakopoulos, R. (2016).
\newblock Network-flow based algorithms for scheduling production in
  multi-processor open-pit mines accounting for metal uncertainty.
\newblock {\em European Journal of Operational Research}, 250(1):273--290.

\bibitem[Lamghari and Dimitrakopoulos, 2020]{LAMGHARI2020104590}
Lamghari, A. and Dimitrakopoulos, R. (2020).
\newblock Hyper-heuristic approaches for strategic mine planning under
  uncertainty.
\newblock {\em Computers and Operations Research}, 115:104590.

\bibitem[Lamghari and Dimitrakopoulos, 2022]{lamghari2022adaptive}
Lamghari, A. and Dimitrakopoulos, R. (2022).
\newblock An adaptive large neighborhood search heuristic to optimize mineral
  value chains under metal and material type uncertainty.
\newblock {\em International Journal of Mining, Reclamation and Environment},
  36(1):1--25.

\bibitem[Lamghari et~al., 2015]{lam15}
Lamghari, A., Dimitrakopoulos, R., and Ferland, J.~A. (2015).
\newblock A hybrid method based on linear programming and variable neighborhood
  descent for scheduling production in open-pit mines.
\newblock {\em Journal of Global Optimization}, 63:555--582.

\bibitem[Newman et~al., 2010]{new10}
Newman, A.~M., Rubio, E., Caro, R., Weintraub, A., and Eurek, K. (2010).
\newblock A {R}eview of {O}perations {R}esearch in {M}ine {P}lanning.
\newblock {\em Interfaces}, 40:222--245.

\bibitem[Osanloo et~al., 2008]{osa08}
Osanloo, M., Gholamnejad, J., and Karimi, B. (2008).
\newblock Long-term open pit mine production planning: a review of models and
  algorithms.
\newblock {\em International Journal of Mining, Reclamation, and Environment},
  22:3--35.

\bibitem[O’Sullivan and Newman, 2015]{o2015optimization}
O’Sullivan, D. and Newman, A. (2015).
\newblock Optimization-based heuristics for underground mine scheduling.
\newblock {\em European Journal of Operational Research}, 241(1):248--259.

\bibitem[Sari and Kumral, 2016]{Sari2016-bc}
Sari, Y.~A. and Kumral, M. (2016).
\newblock An improved meta-heuristic approach to extraction sequencing and
  block routing.
\newblock {\em J. South Afr. Inst. Min. Metall.}, 116(7):673--680.

\bibitem[Sen\'{e}cal and Dimitrakopoulos, 2020]{senecal2020}
Sen\'{e}cal, R. and Dimitrakopoulos, R. (2020).
\newblock Long-term mine production scheduling with multiple processing
  destinations under mineral supply uncertainty, based on multi-neighbourhood
  tabu search.
\newblock {\em International Journal of Mining, Reclamation and Environment},
  34(7):459--475.

\bibitem[Smith and Wicks, 2014]{doi:10.1287/inte.2014.0737}
Smith, M.~L. and Wicks, S.~J. (2014).
\newblock Medium-term production scheduling of the lumwana mining complex.
\newblock {\em Interfaces}, 44(2):176--194.

\bibitem[Zeng et~al., 2021]{ZENG2021102274}
Zeng, L., Liu, S.~Q., Kozan, E., Corry, P., and Masoud, M. (2021).
\newblock A comprehensive interdisciplinary review of mine supply chain
  management.
\newblock {\em Resources Policy}, 74:102274.

\end{thebibliography}

 \APPENDIX

\section{Model Formulation}\label{sec:Model}

A flow network is defined for each pit $m \in \mathcal{M}$. Each network contains a source node,  also denoted by $m$, representing the pit from which blocks are extracted, and a number of destination nodes. These destinations, capturing stockpiles, crushers, processing plants, and dumps, are connected to a virtual node $r$, for each saleable product $r \in \mathcal{R}$. Tables \ref{tab:SetsIndices}--\ref{tab:Parameters} identify and define the sets, variables, and parameters involved in our  model.

\begin{table}[tbp]
\centering
\TABLE{Sets and indices involved in our long-term mine planning model.\label{tab:SetsIndices}}
{\begin{tabular}{lp{350pt}}
\toprule
\textbf{Set} & \textbf{Description} \\
\toprule
    $\mathcal{M}$ &   Pits, indexed by $m$. \\

    $\mathcal{R}$ &  Saleable products produced, indexed by $r$.\\
  
    $\mathcal{N}^{prod}$ &  Product nodes in the overall flow network. \\
   
    $\mathcal{F}_m$ &  Flow network for pit $m \in \mathcal{M}$ with nodes $\mathcal{N}_m$  and arcs  $\mathcal{A}_{m}$. A directed arc $(i,j) \in \mathcal{A}_m$ connects nodes $i \in \mathcal{N}_m$ and $j \in \mathcal{N}_m \cup \mathcal{N}^{prod}$. The arc $(m,j) \in \mathcal{A}_m$ denotes the path from the blocks in pit $m$ to  destination $j$. A single source node is defined to represent the set of blocks in a pit. \\
   
    $\mathcal{N}^{stock}_m$ &  Nodes in the flow network for pit $m \in \mathcal{M}$ that represent stockpiles. \\
  
    $\mathcal{B}$ &  Blocks in the model, indexed by $b$.\\
  
    $\mathcal{B}_m$ &  Blocks in pit $m \in \mathcal{M}$, indexed by $b$.\\
 
    $\mathcal{D}$ &  Mining precedences between blocks in $\mathcal{B}$, where each $(i, j) \in \mathcal{D}$ indicates that the extraction of block $i \in \mathcal{B}$ cannot begin until block $j \in \mathcal{B}$ has been depleted.\\
    
    $\mathcal{L}_{b}$ &  Material types in block $b \in \mathcal{B}_m$ of pit $m \in \mathcal{M}$, indexed by $l$. Where a variable, set or parameter relates to a block $b \in \mathcal{B}_m$, we do not include the index $m$ in its subscript, but rather let $b$ imply the pit.\\
    
    $\mathcal{P}_{bl}$ & Parcels of type $l \in \mathcal{L}_{b}$ in block $b \in \mathcal{B}_m$, pit $m \in \mathcal{M}$, indexed by $p$. \\
    
    $\mathcal{T}$ & Time periods in the planning horizon, indexed by $t$.\\
   
    $\mathcal{J}$ & A collection of \textit{sets of pits}, indexed by $k$, where for each pit set $\mathcal{J}_k$ a minimum production constraint is present. \\
    
    $\mathcal{E}$ & Mineral elements of interest, indexed by $e$.\\
\bottomrule
\end{tabular}}
{}
\end{table}

\begin{table}[tbp]
\centering
\TABLE{Variables involved in our long-term mine planning model.\label{tab:Variables}}
{\begin{tabular}{lp{350pt}}
\toprule
\textbf{Variable} & \textbf{Description} \\
\toprule
    $f_{pijt}$ & Fraction of parcel $p \in \mathcal{P}_{bl}$, of type $l \in \mathcal{L}_b$, from block $b \in \mathcal{B}_m$ in pit $m \in \mathcal{M}$,  that flows along the path $(i,j) \in \mathcal{A}_m$ in period $t \in \mathcal{T}$.  Where a variable relates to a parcel $p$, we exclude indices $b$, $l$, and $m$ from its subscript, letting $p$ imply the block, type, and pit. \\

    $x_{bt}$ & Fraction of block $b \in \mathcal{B}_m$, pit $m \in \mathcal{M}$, extracted by the end of $t \in \mathcal{T}$.\\

    $y_{bt}$ & Binary  that takes on a value of 1  if and only if block $b \in \mathcal{B}_m$, pit $m \in \mathcal{M}$, has been completely extracted by the end of  $t \in \mathcal{T}$.\\

    $z_{bt}$ & Binary  that takes on a value of 1 if and only if extraction of block $b \in \mathcal{B}_m$, pit $m \in \mathcal{M}$, has commenced by the end of  $t \in \mathcal{T}$. \\

    $s_{ipt}$ & Fraction of parcel $p \in \mathcal{P}_{bl}$ of type $l \in \mathcal{L}_b$, from block $b \in \mathcal{B}_m$, pit $m \in \mathcal{M}$, in stockpile $i \in \mathcal{N}^s_m$ in period $t \in \mathcal{T}$.\\

    $w_{mt}$ & Binary that takes on a value of 1 if and only if pit $m \in \mathcal{M}$ is opened in period $t \in \mathcal{T}$. \\

    $v_{mt}$ & Binary that takes on a value of 1 in period $t \in \mathcal{T}$ if and only if pit $m \in \mathcal{M}$ has been opened in, or prior to, $t$.  \\
\bottomrule
\end{tabular}}
{}
\end{table}

\begin{table}[tbp]
\centering
\TABLE{Parameters involved in our long-term mine planning model.\label{tab:Parameters}}
{\begin{tabular}{lp{350pt}}
\toprule
\textbf{Parameter} & \textbf{Description} \\
\toprule
$\pi_t$ & Discount factor for period $t \in \mathcal{T}$.\\

$\tau_p$ & Tonnage of parcel $p$.\\

$\sigma_r$ & Revenue per ton of product $r \in \mathcal{R}$.\\

$\alpha_p$ & Extraction cost per ton for parcel $p$.\\

$\alpha^{cap}_m$ & Capex investment cost associated with opening pit $m \in \mathcal{M}$.\\

$\beta_{ep}$ & Concentration (level) of mineral element $e \in \mathcal{E}$ in parcel $p$.\\[4pt]

$\bar{\rho}_{er}$ &  Maximum allowed level of mineral element $e \in \mathcal{E}$ in product $r \in \mathcal{R}.$\\[4pt]

$\underline{\rho}_{er}$ & Minimum required level of mineral element $e \in \mathcal{E}$ in product $r \in \mathcal{R}$.\\[4pt]

$\bar{\lambda}_n$ & Maximum tonnage permitted to exit node $n$, where $n$ is a node in the overall network flow. Where $n$ is the source node of a pit's flow network, $\bar{\lambda}_n$ is the maximum mining capacity for that pit. \\[4pt]

$\underline{\lambda}_{kt}$ & Minimum tonnage to be mined from the pits $k \in \mathcal{J}$ in period $t \in \mathcal{T}$.\\[4pt]

$\bar{\gamma}_{it}$ & Maximum capacity of stockpile $i \in \mathcal{N}^{stock}_m$, pit $m \in \mathcal{M}$, in $t \in \mathcal{T}$. \\
\bottomrule
\end{tabular}}
{}
\end{table}

\subsection{Objective}

Our objective is to maximize the net present value (NPV) of the mining complex over a long-term planning horizon,    Equation \eqref{eqn:Objective}, where $\mathsf{Rev}_{mt}$ denotes the revenue associated with products formed from ore in pit $m \in \mathcal{M}$ in period $t$ (Equation \ref{eqn:Revenue}), and $\mathsf{Cost}_{mt}$ the total costs incurred from the extraction and processing of material from pit $m$ (Equation \ref{eqn:Cost}). The revenue generated from pit $m \in \mathcal{M}$ in period $t \in \mathcal{T}$ considers the contribution from each parcel  $p \in \mathcal{P}_{bl}$, of each type $l \in \mathcal{L}_b$, from each block in the pit $b \in \mathcal{B}_m$, to each saleable product $r \in \mathcal{R}$, and the revenue generated per ton of that product. The cost associated with  pit $m$ in period $t$ includes capex costs incurred when opening the pit, and extraction costs. The former is incurred in the period in which the pit is opened (where $w_{mt} = 1$).  Extraction costs consider the material extracted from the pit, moving along paths from its source node $m$ to destinations in its flow network.
\begin{linenomath*}
\begin{equation}
    \max \sum_{t \in \mathcal{T}} \pi_t \left( \sum_{m \in \mathcal{M}} \mathsf{Rev}_{mt} -  \mathsf{Cost}_{mt}\right)
    \label{eqn:Objective}
\end{equation}
\end{linenomath*}
where: 
\begin{linenomath*}
\begin{multline}
\mathsf{Rev}_{mt} = \sum_{b \in \mathcal{B}_m} 
\sum_{l \in \mathcal{L}_{b}} \sum_{p \in \mathcal{P}_{bl}} 
\sum_{\substack{(i,r) \in \mathcal{A}_{m} \\ \mid r \in \mathcal{N}^{prod}}} \tau_{p} \, 
\sigma_{r} \, f_{pirt} \label{eqn:Revenue} \\
\forall t \in \mathcal{T}, \; \forall m \in \mathcal{M} \nonumber \\
\mathsf{Cost}_{mt} =  \alpha^{cap}_m w_{mt}  + \sum_{b \in \mathcal{B}_m} 
\sum_{l \in \mathcal{L}_{b}} \sum_{p \in \mathcal{P}_{bl}} 
\sum_{(m,j) \in \mathcal{A}_{m}} \tau_{p}  
\alpha_p  f_{pmjt} \label{eqn:Cost} \\
\forall t \in \mathcal{T},  \forall m \in \mathcal{M} \nonumber
\end{multline} 
\normalsize
\end{linenomath*}

\subsection{Constraints}\label{sec:ModelConstraints}
\begin{linenomath*}
Blending constraints \eqref{c:proddefLB}--\eqref{c:proddefUB} enforce lower and upper bounds on the concentration of each mineral element  $e \in \mathcal{E}$ in each formed product $r \in \mathcal{R}$. Constraint set \eqref{c:maxextractnode} enforces capacities on the tonnage of material exiting any given intermediate node $n$ in our overall flow network. This excludes source nodes in the flow network for each pit, and product nodes. 
Constraint set \eqref{c:maxextractmine} enforces mining capacities in each pit $m \in \mathcal{M}$. Recall that $v_{mt}$ is a binary indicating whether pit $m$ has been opened by, or in, period $t \in \mathcal{T}$. When set to 1, the pit's mining capacity becomes $\bar{\lambda}_{m}$.  Constraint set \eqref{c:minprod} enforces minimum production constraints across collections of pits, $k \in \mathcal{J}$. 
\begin{multline} 
    \sum_{m \in \mathcal{M}}
    \sum_{b \in \mathcal{B}_m} 
    \sum_{l \in \mathcal{L}_{b}}     
    \sum_{p \in \mathcal{P}_{bl}} 
    \sum_{(i,r) \in \mathcal{A}_{m}} 
    { (\beta_{ep} - \underline{\rho}_{er}) \, \tau_{p} \, f_{pirt}} \geq 0 \label{c:proddefLB}\\
    \; \forall r \in \mathcal{N}^{prod}, \; \forall e \in \mathcal{E}, \; \forall t \in \mathcal{T}\nonumber\\ 
    \sum_{m \in \mathcal{M}}
    \sum_{b \in \mathcal{B}_m} 
    \sum_{l \in \mathcal{L}_{b}}     
    \sum_{p \in \mathcal{P}_{bl}} 
    \sum_{(i,r) \in \mathcal{A}_{m}} 
    { (\beta_{ep} - \bar{\rho}_{er}) \, \tau_{p} \, f_{pirt}} \leq 0 \label{c:proddefUB}\\
    \; \forall r \in \mathcal{N}^{prod}, \; \forall e \in \mathcal{E}, \; \forall t \in \mathcal{T}\nonumber\\
    \sum_{b \in \mathcal{B}_m} 
    \sum_{l \in \mathcal{L}_{b}}     
    \sum_{p \in \mathcal{P}_{bl}} 
    \sum_{(n,j) \in \mathcal{A}_{m} } 
     \tau_{p} \, f_{pnjt} \leq \bar{\lambda}_n \label{c:maxextractnode}\\ 
     \forall m \in \mathcal{M}, \;
     \forall n \in \mathcal{N}_{m} \mid n \neq m, \;
     \forall t \in \mathcal{T}\nonumber\\
    \sum_{b \in \mathcal{B}_m} 
    \sum_{l \in \mathcal{L}_{b}}     
    \sum_{p \in \mathcal{P}_{bl}} 
    \sum_{(m,j) \in \mathcal{A}_{m}} 
     \tau_{p} \, f_{pmjt} \leq \bar{\lambda}_{m} \, v_{mt} \label{c:maxextractmine}\\
     \forall m \in \mathcal{M}, \;
     \forall t \in \mathcal{T} \nonumber\\ 
    \sum_{m \in \mathcal{J}_k}
    \sum_{b \in \mathcal{B}_m} 
    \sum_{l \in \mathcal{L}_{b}}     
    \sum_{p \in \mathcal{P}_{bl}} 
    \sum_{(m,j) \in \mathcal{A}_{m}} 
     \tau_{p} \, f_{pmjt} \geq  \underline{\lambda}_{kt} \label{c:minprod} \\
     \forall k \in \mathcal{J},  \;
     \forall t \in \mathcal{T} \nonumber
 \end{multline}


Recall that $w_{mt}$ is a binary variable signalling whether the pit $m \in \mathcal{M}$ was opened in period $t$. This binary variable triggers the 
binary $v_{mt}$, which in turn activates extraction capacity for pit $m$ from period $t$ onward.
\begin{multline}
 v_{mt} \leq \sum_{\tau=1,2, \ldots t} w_{m\tau}  \\
 \forall m  \in \mathcal{M}, \; \forall t \in \mathcal{T} \nonumber
\end{multline}

Mass balance constraints \eqref{c:mbalance} are defined for all intermediate nodes in our overall flow network, excluding source nodes and those representing stockpiles (nodes $\mathcal{N}^{stock}_{m}$ for pit $m \in \mathcal{M}$). 
Mass balance constraints \eqref{c:spmbalance} are also defined for each node in our flow network representing a stockpile to keep track of the fraction of each parcel present in the stockpile in any given time period. Constraint set \eqref{c:splink} ensures that the fraction of a parcel leaving a stockpile in period $t \in \mathcal{T}$ does not exceed the fraction of the parcel that was present in the stockpile in period $t - 1$.  Stockpile capacities are enforced by constraint set \eqref{c:spbnd}, where $\bar{\gamma}_{it}$ denotes the capacity of stockpile $i \in \mathcal{N}^{stock}_m$ in period $t$.
\begin{multline} 
   \sum_{(i,n) \in \mathcal{A}_{m} } 
    f_{pint} = 
    \sum_{(n,j) \in \mathcal{A}_{m} } {  f_{pnjt} }  \label{c:mbalance} \\
    \forall t \in \mathcal{T}, 
    \forall m \in \mathcal{M}, 
    \forall b \in \mathcal{B}_m,     
    \forall l \in \mathcal{L}_{b}, \nonumber \\
    \forall p \in \mathcal{P}_{bl}, 
    \forall n \in \mathcal{N}_{m} \setminus \mathcal{N}^{stock}_{m} \cup \{m\}  \nonumber \\
    s_{ipt} = s_{ip,t-1} +  
    \sum_{(j,i) \in \mathcal{A}_{m}} f_{pjit} -
       \sum_{(i,j) \in \mathcal{A}_{m}}{f_{pijt}}  \label{c:spmbalance} \\ 
     \forall m \in \mathcal{M},  
      \forall i \in \mathcal{N}^{stock}_m, \;
      \forall b \in \mathcal{B}_m, \forall l \in \mathcal{L}_{b}, \nonumber \\
      \forall p \in \mathcal{P}_{bl}, 
      \forall t \in \mathcal{T} \mid s_{ip0}= 0 \nonumber \\
    \sum_{(i,j) \in \mathcal{A}_{m}}  f_{pijt} \leq s_{ip,t-1}    \label{c:splink}\\ 
    \forall m \in \mathcal{M}, \;
    \forall i \in \mathcal{N}^{stock}_m, \;
    \forall b \in \mathcal{B}_m, \;     
    \forall l \in \mathcal{L}_{b}, \nonumber  \\
      \forall p \in \mathcal{P}_{bl}, 
    \forall t \in \mathcal{T} \mid t \geq 2 \nonumber\\
    \sum_{b \in \mathcal{B}_m}{ \sum_{l \in \mathcal{L}_{b}} 
    { \sum_{p \in \mathcal{P}_{bl}}} \tau_{p} s_{ipt}} \leq  \bar{\gamma}_{it}  \label{c:spbnd}\\
    \forall m \in \mathcal{M}, \;
    \forall i \in \mathcal{N}^{stock}_m, \;
    \forall t \in \mathcal{T}  \nonumber
\end{multline}

Constraint set \eqref{eq:block1} enforces mining precedences between blocks. For a given precedence relationship $(i, j) \in \mathcal{D}_m$ between blocks $i$ and $j$ in pit $m$, extraction of block $i$ cannot begin until block $j$ has been depleted. We cannot progress the mining of a block $b \in \mathcal{B}_m$ until extraction of the block has begun, $z_{bt} = 1$ \eqref{eq:block2}. A block $b \in \mathcal{B}_m$ is not depleted until the fraction of the block mined by a time period, $x_{bt}$, is equal to 1 \eqref{eq:block3}. Constraint sets \eqref{c:coherence1}--\eqref{c:coherence3} enforce coherence in the state of a block  over time. Constraint set \eqref{eq:linkblocks} connects the material flow variables for a block to its progression variables. Constraint set \eqref{eq:parceltons} restricts the mining of a parcel to its available tonnage. 
\begin{multline} 
    z_{it} \leq y_{jt}  \quad
    \forall \left(i, j\right) \in \mathcal{D}, \; \forall t \in \mathcal{T} \label{eq:block1}\\
    x_{bt} \leq z_{bt} \quad 
    \forall b \in \mathcal{B}, \; \forall t \in \mathcal{T}\label{eq:block2}\\
    y_{bt} \leq x_{bt} \quad 
    \forall b \in \mathcal{B}, \; \forall t \in \mathcal{T}\label{eq:block3}\\ 
    x_{bt} \leq x_{b,t+1} \quad  
     \forall b \in \mathcal{B}, \; 
    \forall t \in \mathcal{T}  \label{c:coherence1}\\
    y_{bt} \leq y_{b,t+1} \quad  
     \forall b \in \mathcal{B}, \; 
    \forall t \in \mathcal{T} \\
    z_{bt} \leq z_{b,t+1} \quad  
    \forall b \in \mathcal{B}, \; 
    \forall t \in \mathcal{T} \label{c:coherence3}\\
     \sum_{(m,j) \in \mathcal{A}_{m}}  \sum_{l \in \mathcal{L}_{b}} {\sum_{p \in \mathcal{P}_{bl}}f_{pmjt}}   =  
    x_{bt} - x_{b,t-1}  \label{eq:linkblocks} \\
    \forall m \in \mathcal{M}, \; 
    \forall b \in \mathcal{B}_m, \; 
    \forall t \in \mathcal{T} \mid x_{b0} = 0 \nonumber \\
    \sum_{t \in \mathcal{T}} \sum_{(m,j) \in \mathcal{A}_{m}}  f_{pmjt}  \leq \tau_p  \label{eq:parceltons}\\
    \forall m \in \mathcal{M}, \; 
    \forall b \in \mathcal{B}_m, \; 
    \forall l \in \mathcal{L}_b, \;
    \forall p \in \mathcal{P}_{bl} \nonumber
\end{multline}

\end{linenomath*}

\section{Detailed Results}\label{sec:DetailedResults}

We evaluate the use of sliding windows for generating initial feasible solutions across varying window size ($W$), relaxed horizon length ($H$), and overlap ($O$). Tables \ref{tab:OTSlidingWindows} and \ref{tab:PilbaraSlidingWindows} report the quality of solutions found by sliding windows, in terms of their gap to either a known optimal solution or a best known upper bound, for the Oyu Tolgoi and Pilbara models. We vary $W$ from 1 to 7, $H$ from 0 to 4, and $O$ from 0 to 2. A `--' indicates either that the combination of parameter values could not be used in conjunction, or that sliding windows did not yield a feasible solution.

\begin{table}[t]
\centering
\small
\TABLE{\textbf{Oyu Tolgoi}: Quality of solutions formed by sliding windows, recording the gap \% to the best bound found after 148 hours of solving the full model, and solve time (s)  for varying $W$, $H$, and $O$. MIPs are terminated at a MIP gap of 0.1\%. \label{tab:OTSlidingWindows}}
{\begin{tabular}{c|c|c|c|c|c|c|c|c|c|c}
\toprule
  & \multicolumn{2}{c|}{$H = 0, \,O = 0$} &  \multicolumn{2}{c|}{$H = 0, \,O = 1$} &  \multicolumn{2}{c|}{$H = 0, \,O = 2$} &  \multicolumn{2}{c|}{$H = 2, \,O = 1$} &  \multicolumn{2}{c}{$H = 4, \,O = 1$}\\
 \cline{2-11}
$\,W\,$   & Gap  (\%)  &  Time (s)  & Gap  (\%)  &  Time (s)  & Gap (\%)   &  Time (s) & Gap (\%)   &  Time (s)& Gap (\%)   &  Time (s)\\ 
\midrule
1  & 99.9 & 152 & -- & -- & -- & -- & -- & -- & -- & -- \\
2  & 99.9 & 110 & 99.9 & 150 & -- & -- & 99.9 & 162 & 37.7 & 168\\
3  & 99.9 & 97 & 99.9 & 112 & 99.9 & 144 & 99.9 & 122 & 37.8 & 125\\
4  & 99.9 & 92 & 99.9 & 98 & 99.9 & 103 & 48.3 & 105 & 27.1 & 108\\
5  & 68.3 & 85 & 59.6 & 92 & 48.4 & 90 & 27.1 & 98 & 27.1 & 102\\
6  & 27.1 & 87 & 27.1 & 85 & 27.1 & 84 & 27.1 & 89 & 27.1 & 88\\
7  & 27.1 & 78 & 27.1 & 87 & 27.1 & 87 & 27.1 & 90 & 27.1 & 90\\
\bottomrule
\end{tabular}}
{}
\end{table}

\begin{table}[t]
\centering
\TABLE{\textbf{Pilbara}: Quality of solutions formed by sliding windows, recording the gap \% to the best bound found after 148 hours of solving the full model, and solve time (min)  for varying $W$, $H$, and $O$. MIPs are  terminated at a MIP gap of 10\%, or after 5000s.
\label{tab:PilbaraSlidingWindows}}
{\begin{tabular}{c|c|c|c|c|c|c|c|c|c|c}
\toprule
  & \multicolumn{2}{c|}{$H = 0, \,O = 0$} &  \multicolumn{2}{c|}{$H = 0, \,O = 1$} &  \multicolumn{2}{c|}{$H = 0, \,O = 2$} &  \multicolumn{2}{c|}{$H = 2, \,O = 1$} &  \multicolumn{2}{c}{$H = 4, \,O = 1$}\\
 \cline{2-11}
$W$   & Gap  (\%)  &  Time   & Gap  (\%)  &  Time  & Gap (\%)   &  Time  & Gap (\%)   &  Time & Gap (\%)   &  Time \\  
   &   &  (min)  &  &  (min)  &    &   (min) &   &   (min)&   &  (min)\\ 
\midrule
1  & 12.2 & 31.7 & -- & -- & -- & -- & -- & -- & -- & --\\
2  & 9.9 & 35 & 8.9 & 51.7 & -- & -- & 9.6 & 260 & -- & --\\
3  & 9.3 & 45 & 8.9 & 50 & 8.2 & 120 & 8.3 & 238.3 & -- & --\\
4  & 9.6 & 68.3 & 8.2 & 81.7 & 9.6 & 121.7 & -- & -- & -- & --\\
\bottomrule
\end{tabular}}
{}
\end{table}

\begin{table}[!tb]
\centering
\TABLE{\textbf{Oyu Tolgoi}: Average LNS run time (min), and quality of solutions found (\% gap to best known upper bound) across 10 runs. Initial solution from sliding windows has a gap of 27.1\%. Solutions compared given varying focal block selection (MD,  MIX, OBJ, and RAND), neighbourhood size (N), variable fixing  (S/D and S/D + F), and use of RINS (Yes/No).  LNS MIPs are terminated at a gap of 0.05\%, or after 600s. LNS is terminated after 12 hrs or when the improvement rate fell below 1\%.  Best settings in bold. Baseline setting: 2.
\label{tab:OTNSizeVary}}
{\begin{tabular}{l|c||c|c|c|c|c}
\toprule
$\#$ &   $N$    & Focal   &  Use of      & Fixing  & Avg Time    & Avg Gap  \\
 &     &   Block   &   RINS     &  Method  & (Min, Max) (min)   & (Min, Max)  (\%) \\
 \midrule
 \multicolumn{7}{l}{Neighbourhood formation strategy: \textit{Timing} and \textit{Blending}} \\
  \midrule
1 & 50 & OBJ & Yes & S/D & 38.3 (25,50)   & 17.5 (12.5,26.0)    \\
\bf 2$^*$ & \bf 100 & \bf OBJ & \bf Yes & \bf S/D & \bf 258.3, (161.7,315)   &   \bf 10.3 (9.1,11.3) \\
3 &  100 &  OBJ &  Yes &  S/D + F &  251.7, (186.7,338.3)   &    10.7 (8.9,11.3) \\
4 & 100 & OBJ & No & S/D & 205, (91.7,318.3)   &   13.0 (9.1,24.0)\\
5 & 150 & OBJ & Yes & S/D & 378.3, (343.3,461.7)   &   10.5 (9.0,11.5) \\
\midrule
\bf 6 & \bf 100 & \bf MD & \bf Yes & \bf S/D & \bf 323.3, (290,353.3)   &   \bf 9.9 (7.8,11.2)\\
\bf 7 & \bf 100 & \bf RAND & \bf Yes & \bf S/D &  \bf 338.3, (300,391.7)   &   \bf 9.6 (7.9,10.8)\\
\bf 8 & \bf 100 & \bf MIX & \bf Yes & \bf S/D & \bf 315, (280,353.3)   &  \bf 10.0 (8.4,10.7) \\
\midrule
\midrule
 \multicolumn{7}{l}{Neighbourhood formation strategy: \textit{Timing} only} \\
  \midrule
9 &   100 & OBJ & Yes & S/D &   133.3, (83.3,291.7)   &   11.4 (10.2,13.1) \\
 \midrule
 \midrule
 \multicolumn{7}{l}{Neighbourhood formation strategy: \textit{Blending} only} \\
  \midrule 
10 &     100 & OBJ & Yes & S/D &   285, (211.7,355)   &   10.5 (8.8,12.0)  \\
    \midrule 
     \multicolumn{7}{l}{Neighbourhood formation strategy: None} \\
  \midrule 
11 &    100 & OBJ & Yes & S/D & 93.3, (28.3,176.7)   &   13.9 (11.1,26.0)\\
 \bottomrule
\end{tabular}}
{}
\end{table}

For Oyu Tolgoi, Table \ref{tab:OTNSizeVary} reports the average quality of solutions found by LNS  across varying neighbourhood sizes, focal block selection methods, variable fixing strategies, and neighbourhood formation strategies. We report solution quality in terms of the gap to the best known upper bound on the objective, found after solving the full model with Gurobi for 148 hours, warm started with the same initial solution used as  the starting point for LNS.  The unfix window $UW$ is set to the length of the planning horizon. The entries in Table \ref{tab:OTNSizeVary} start with an initial solution with a gap of 27.1\%, obtained by sliding windows using $H = 2,$ $O = 1,$ and $W = 5$. LNS MIPs are solved to a gap of 0.05\%, or to a time limit of 600s. LNS is terminated after either 12 hours or when the improvement rate falls below 1\%.

Table \ref{tab:PilbaraLNS} reports the average quality of solutions found by LNS to the Pilbara model across varying parameter settings. Using a fixed neighbourhood size of 200 blocks, we report the average quality of solutions found by LNS, for each setting, in terms of their gap to the best known upper bound on the optimal objective, over 10 differently seeded runs, each utilising 10 threads of computation. LNS MIPs were solved to a gap of 0.05\%, or a time limit of 400s. Each LNS run was terminated after 12 hours. The entries in Table \ref{tab:PilbaraLNS} have started with an initial solution with a gap of 8.2\%, found by sliding windows with $H = 0$, $O = 1$ and $W = 4$. The best settings are in bold.

 \begin{table}[tpb]
\centering
\TABLE{\textbf{Pilbara}: Average quality of LNS solutions over 10 runs (\% gap to best upper bound found after 148 hours of solving the full model). Columns are defined as per Table \ref{tab:OTNSizeVary}, with $UW$ denoting the size of the unfix window used for each MIP solve. Initial solution with a gap of 8.2\% found by sliding windows ($H = 0,$ $O = 1,$ and $W = 4$). Neighbourhood size of 200 blocks.  LNS MIPs terminated at a gap of 0.05\%, or after 400s. LNS is terminated after 12 hrs. Best settings in bold. Baseline setting: 8.
\label{tab:PilbaraLNS}}
{\begin{tabular}{l|l|c|c|c|c|c}
\toprule
&  Neighbourhood    &  Focal Block & Fixing &  Use of  & Unfix &  Avg Gap   \\
$\#$ &  Strategies Used  &  Selection & Method &      RINS & Window ($UW$)  &    (Min, Max)  (\%)      \\
 \midrule
1 & None              & OBJ &  S/D+F &  No & 2 & 4.7 (4.5, 4.9) \\
\midrule
2 & All & OBJ & S/D & No & 2 &  8.2 (8.2,8.2) \\
3 & All & OBJ & S/D + F & Yes & 2 &  4.8 (4.5,5.1) \\
4 & All & OBJ & S/D + F & No & 4  &4.9 (4.6, 5.3) \\
\midrule
\bf 5 & \bf Trigger           & \bf OBJ & \bf  S/D+F &  \bf No &  \bf 2 & \bf 4.5 (4.4, 4.6) \\ 
6 & Timing, Blending  & OBJ &  S/D+F &  No &  2 & 4.7 (4.6,4.9)  \\ 
7 & Pit Links, Timing & OBJ &  S/D+F &  No &  2 & 4.6 (4.3, 4.9) \\ 
\midrule
\bf 8$^*$ & \bf All               & \bf OBJ & \bf S/D+F & \bf No & \bf 2 & \bf 4.5 (4.2, 4.7) \\ 
\bf 9 & \bf $-$ Blending      & \bf OBJ & \bf S/D+F & \bf No & \bf 2 & \bf 4.5 (4.2,4.8)\\ 
\bf 10 & \bf $-$ Pit Links     & \bf OBJ & \bf S/D+F & \bf No & \bf 2 & \bf 4.5 (4.3,4.8)\\ 
\midrule
\bf 11 & \bf All               &  \bf MD & \bf S/D+F & \bf No & \bf 2 & \bf 4.5 (4.3, 4.8) \\ 
\bf 12 & \bf All               & \bf RAND & \bf S/D+F & \bf No & \bf 2 & \bf 4.5 (4.2, 4.7)\\ 
\bf 13 & \bf All               & \bf MIX  & \bf S/D+F &\bf No & \bf 2 & \bf4.5 (4.3, 4.8)\\ 
 \bottomrule
\end{tabular}}
{}
\end{table}

\end{document}